\documentclass{article}
\usepackage{CJK,CJKnumb,CJKulem,times,dsfont,ifthen,mathrsfs,latexsym,amsfonts,color}
\usepackage{amsmath,amsthm,makeidx,fontenc,amssymb,bm,graphicx,psfrag,listings,curves,extarrows,enumitem}

\usepackage[title]{appendix}

\usepackage[colorlinks=true,urlcolor=black,
citecolor=black,linkcolor=black,linktocpage,pdfpagelabels,
bookmarksnumbered,bookmarksopen]{hyperref}
\usepackage[colorinlistoftodos]{todonotes}

\usepackage{geometry}
\usepackage{indentfirst}
\usepackage{amssymb}
\makeatletter

\newcommand{\Rmnum}[1]{\expandafter\@slowromancap\romannumeral #1@}
\makeatother

\newtheorem{theorem}{Theorem}[section]
\newtheorem{definition}[theorem]{Definition}
\newtheorem{lemma}[theorem]{Lemma}
\newtheorem{corollary}[theorem]{Corollary}
\newtheorem{proposition}[theorem]{Proposition}
\newtheorem{remark}[theorem]{Remark}

\newcommand{\al}{\alpha}
\newcommand{\ga}{\gamma}
\newcommand{\Ga}{\Gamma}
\newcommand{\dl}{\Delta}
\newcommand{\e}{\varepsilon}

\newcommand{\iy}{\infty}
\newcommand{\q}{\theta}

\newcommand{\la}{\lambda}
\newcommand{\vp}{\varphi}
\newcommand{\pa}{\partial}
\newcommand{\ti}{\tilde}

\newcommand{\ra}{\rightarrow}
\newcommand{\rh}{\rightharpoonup}

\newcommand{\lab}{\label}
\newcommand{\f}{\frac}
\newcommand{\bt}{\begin{theorem}}
\newcommand{\et}{\end{theorem}}
\newcommand{\bl}{\begin{lemma}}
\newcommand{\el}{\end{lemma}}
\newcommand{\bd}{\begin{definition}}
\newcommand{\ed}{\end{definition}}
\newcommand{\bc}{\begin{corollary}}
\newcommand{\ec}{\end{corollary}}
\newcommand{\bp}{\begin{proof}}
\newcommand{\ep}{\end{proof}}
\newcommand{\bx}{\begin{example}}
\newcommand{\ex}{\end{example}}
\newcommand{\bi}{\begin{exercise}}
\newcommand{\ei}{\end{exercise}}
\newcommand{\br}{\begin{remark}}
\newcommand{\er}{\end{remark}}
\newcommand{\be}{\begin{equation}}
\newcommand{\ee}{\end{equation}}
\newcommand{\bal}{\begin{align}}
\newcommand{\bn}{\begin{enumerate}}
\newcommand{\en}{\end{enumerate}}
\newcommand{\ba}{\begin{align}}
\newcommand{\ea}{\begin{align}}
\newcommand{\bg}{\begin{align*}}
\newcommand{\eg}{\end{align*}}
\newcommand{\bcs}{\begin{cases}}
\newcommand{\ecs}{\end{cases}}

\newcommand{\BR}{{\cal B}}
\newcommand{\CR}{{\cal C}}
\newcommand{\DR}{{\cal D}}

\newcommand{\GR}{{\cal G}}

\newcommand{\LR}{{\cal L}}

\newcommand{\OR}{{\cal O}}
\newcommand{\PR}{{\cal P}} 
\newcommand{\SR}{{\cal S}} 
\newcommand{\TR}{{\cal T}} 
\newcommand{\UR}{{\cal U}}

\newcommand{\C}{{\mathbb C}}
\newcommand{\R}{{\mathbb R}}
\newcommand{\RN}{{\mathbb R^N}}  
\newcommand{\bean}{\begin{eqnarray*}}
\newcommand{\eean}{\end{eqnarray*}}

\newcommand{\supp}{\operatorname{\rm supp}}
\newcommand{\sbr}[1]{\left(#1\right)}
\newcommand{\mbr}[1]{\left[#1\right]}
\newcommand{\lbr}[1]{\left\{#1\right\}}

\newcommand{\s}{\star}

\newcommand{\nm}[1]{\Vert #1 \Vert}
\newcommand{\np}[1]{\Vert #1 \Vert_p}
\newcommand{\nt}[1]{\Vert #1 \Vert_2}
\newcommand{\ns}[1]{\Vert #1 \Vert_{2^*}}
\newcommand{\p}{\ga_{\al+\beta}} 	
\newcommand{\nun}{{\nu_n}}				

\setitemize{itemindent=38pt,leftmargin=0pt,itemsep=-0.4ex,listparindent=26pt,partopsep=0pt,parsep=0.5ex,topsep=-0.25ex}

\numberwithin{equation}{section}

\begin{document}
\theoremstyle{plain}

\title{\bf Existence and asymptotic behavior of normalized ground states for Sobolev critical Schr\"odinger systems}

\date{}
\author{\bf Thomas Bartsch,\; Houwang Li\;\&\;Wenming Zou\thanks{Supported by NSFC(12171265, 11771234)}}

\maketitle
\begin{center}
\begin{minipage}{120mm}
\begin{center}{\bf Abstract }\end{center}
The paper is concerned with the existence and asymptotic properties of normalized ground states of the following nonlinear Schr\"odinger system with critical exponent:
\begin{equation*}
	\left\{\begin{aligned}
	&-\dl u+\la_1 u=|u|^{2^*-2}u+{\nu\al} |u|^{\al-2}|v|^\beta u,\quad \text{in }\RN,\\
	&-\dl v+\la_2 v=|v|^{2^*-2}v+{\nu\beta} |u|^\al |v|^{\beta-2}v,\quad \text{in }\RN,\\
	&\int_\RN u^2=a^2,\;\;\; \int_\RN v^2=b^2,
	\end{aligned}	\right.
\end{equation*}
where $N=3,4$, $\al,\beta>1$, $2<\al+\beta<2^*=\f{2N}{N-2}$. 
We prove that a normalized ground state does not exist for $\nu<0$. 
When $\nu>0$ and $\al+\beta\le2+\f{4}{N}$, we show that the system has a normalized ground state solution for $0<\nu<\nu_0$, 
the constant $\nu_0$ will be explicitly given. 
In the case $\al+\beta>2+\f{4}{N}$ we prove the existence of a threshold $\nu_1\ge 0$ such that a normalized ground state solution exists for $\nu>\nu_1$, and does not exist for $\nu<\nu_1$. We also give conditions for $\nu_1=0$. Finally we obtain the asymptotic behavior of the minimizers as $\nu\to0^+$ or $\nu\to+\iy$.

\vskip0.1in
{\bf Key words:}   Critical Schr\"odinger systems; Normalized solution; Ground states; Asymptotic behaviors.

\vskip0.1in
{\bf 2010 Mathematics Subject Classification:} 35J50, 35J15, 35J60.

\end{minipage}
\end{center}

\vskip0.3in
\section{Introduction and statement of results}
The Schr\"odinger system
\be\lab{201905-1}
\begin{cases}
-i\f{\pa}{\pa t}\Phi_1=\dl \Phi_1+\mu_1 |\Phi_1|^{p-2}\Phi_1+\nu \al|\Phi_1|^{\al-2}|\Phi_2|^{\beta}\Phi_1,\\
-i\f{\pa}{\pa t}\Phi_2=\dl \Phi_2+\mu_2|\Phi_2|^{q-2}\Phi_1+\nu \beta|\Phi_1|^{\al}|\Phi_2|^{\beta-2}\Phi_2,\\
\Phi_j=\Phi_j(x,t)\in \C,\ (x,t)\in \RN\times \R,\ j=1,2,
\end{cases}
\ee
has received a lot of attention in recent years. It appears, for instance, in nonlinear optics, or in mean field models for binary mixtures of Bose-Einstein condensates, or binary gases of fermion atoms in degenerate quantum states (Bose-Fermi mixtures, Fermi-Fermi mixtures). The constants $\mu_1,\mu_2$ and $\nu$ describe the intraspecies and interspecies scattering lengths. The sign of $\nu$ determines whether the interactions are repulsive or attractive, i.e., the interaction is attractive if $\nu>0$, and it is repulsive if $\nu<0$. The most prominent case is that of the coupled Gross-Pitaevskii equations where $p=q=4$ and $\al=\beta=2$. Models for certain ultracold quantum gases however use different exponents; see \cite{Ad, BG-2, BG-3, BG-1, malomed:2008}. 

\vskip0.13in
The ansatz $\Phi_1(x,t)=e^{i\la_1t}u(x)$ and $\Phi_2(x,t)=e^{i\la_2t}v(x)$ for solitary wave solutions of system \eqref{201905-1} leads to the following system of elliptic equations:
\be\lab{equ0}
	\left\{ 	\begin{aligned}
	&-\dl u+\la_1 u = \mu_1 |u|^{p-2}u + \nu\al|u|^{\al-2}|v|^{\beta} u,\quad \text{in }\RN,\\
	&-\dl v+\la_2 v = \mu_2 |v|^{q-2}v + \nu\beta|u|^{\al} |v|^{\beta-2}v,\quad \text{in }\RN,\\
	&u(x),v(x)\to0~~\text{as }|x|\to\iy.
	\end{aligned}	\right.
\ee
Fixing the parameters $\la_1,\la_2>0$ the existence and multiplicity of solutions to \eqref{equ0} has been investigated in the last two decades by many authors. We refer the reader to \cite{Peral etc=CVPDE=2009, Ambrosetti-Colorado=2007,Bartsch-Dancer-Wang:2010, BaWa, Bartsch-Wang-Wei:2007, Chen-Zou=CVPDE=2014,Chen-Zou=TransAMS=2015, Lin-Wei=CMP=2005,Sirakov=CMP=2007} and the references therein.

\vskip0.13in
An important feature of system \eqref{201905-1} is that the $L^2$-norms 
\[
  \int_{\RN}|\Phi_1(t,x)|^2\mathrm{d}x\quad \text{and}\quad\int_{\RN}|\Phi_2(t,x)|^2\mathrm{d}x
\]
of the wave functions are conserved. These have important physical significance: in Bose-Einstein condensates, for example, $|\Phi_1(t,\cdot)|_2$ and $|\Phi_2(t,\cdot)|_2$ represent the number of particles of each component; in the nonlinear optics framework, $|\Phi_1(t,\cdot)|_2$ and $|\Phi_2(t,\cdot)|_2$ represent the power supply. Therefore it is natural to consider \eqref{201905-1} with the constraints 
\[
	\int_\RN u^2=a^2\quad \text{and}\quad \int_\RN v^2=b^2.
\]
Then $\la_1,\la_2$ cannot be prescribed but appear as Lagrange multipliers in a variational approach.
The existence and properties of these {\it normalized solutions} recently has attracted the attention of many reasearchers 
but compared with the results about \eqref{equ0} with $\la_1,\la_2$ prescribed surprisingly little is known. 
This is even true for scalar nonlinear Schr\"odinger equations. 
From a mathematical point of view new difficulties appear in the search of normalized solutions and make this problem mathematically challenging. 
For instance, when one considers the energy functional corresponding to \eqref{equ0} under the mass constraint, in the mass supercritical case $p,q>2+\f4N$ there are bounded Palais-Smale sequences that do not have a convergent subsequence and converge weakly to $0$. We refer the reader to \cite{BaVa, Ikoma-Miyamoto=2020,Ikoma-Tanaka=AdvDE=2019, Jeanjean=1997, Jeanjean-Lu=CVPDE=2020, Jeanjean=2020, Soave=JDE=2020,Soave=JFA=2020, Wei-Wu=2020} for scalar equations, \cite{Bartsch-Jeanjean=2018, Bartsch-Jeanjean-Soave=JMPA=2016, Bartsch-Soave=JFA=2017, Bartsch-Soave=CVPDE=2019, Bartsch-Zhong-Zou=MathAnn=2020, Gou-Jeanjean=2016, Gou-Jeanjean=2018, Li-Zou-1} for systems of two equations, \cite{Mederski-1} for systems of $k$ equations. 

\vskip0.13in
In the present paper we are interested in the existence of ground state normalized solutions of \eqref{equ0} when $p=q=2^*$ is the critical Sobolev exponent,  $2<\al+\beta<2^*$, and $N=3$ or $N=4$. Thus in dimension 3 we deal with a quintic nonlinearity for the self-interaction terms $|u|^{2^*-2}u$ and $|v|^{2^*-2}v$, in dimension 4 with the classical cubic one. We would like to mention that it would be possible to add lower order terms to these self-interactions. In particular, in dimension 3 we could deal with a cubic-quintic nonlinearity which is used in some applications. For simplicity we just deal with homogeneous nonlinearities, and we assume $\mu_1=\mu_2=1$; so we consider the problem 
\be\lab{mainequ}
	\left\{ 	\begin{aligned}
	&-\dl u+\la_1 u=|u|^{2^*-2}u+{\nu\al} |u|^{\al-2}|v|^\beta u,\quad \text{in }\RN,\\
	&-\dl v+\la_2 v=|v|^{2^*-2}v+{\nu\beta} |u|^\al |v|^{\beta-2}v,\quad \text{in }\RN,\\
	&\int_\RN u^2=a^2,\int_\RN v^2=b^2
	\end{aligned}	\right.
\ee
under the assumptions
\be\lab{assumption}
 N\in\{3,4\},\ a,b>0, \text{ and } \al>1,\ \beta>1,\ \al+\beta<2^*.
\ee
It is worth to point out that all the papers about normalized solutions of \eqref{equ0} mentioned above, except \cite{Li-Zou-1}, deal with the Sobolev subcritical case $2<p,q<2^*$. In \cite{Li-Zou-1} only the case $2<p<2^*=q$ has been considered. 
\vskip0.13in
Solutions of \eqref{mainequ} correspond to critical points of the energy functional $I_\nu: H=H^1(\RN)\times H^1(\RN) \to\R$ defined by
\be\lab{energy}
	I_\nu(u,v):=\f{1}{2}\int_\RN \sbr{|\nabla u|^2+|\nabla v|^2}
			-\f{1}{2^*}\int_{\RN} \sbr{|u|^{2^*}+|v|^{2^*}}  -{\nu}\int_{\RN}|u|^\al|v|^\beta
\ee
and constrained to the $L^2$-torus
\[
	\TR(a,b) := \lbr{(u,v)\in  H^1(\RN)\times H^1(\RN):\nt{u}=a,~\nt{v}=b}.
\]
The parameters $\la_1,\la_2\in\R$ will appear as Lagrange multipliers. It is easy to see that $I_\nu$ is of class $C^1$, and that it is unbounded from below on $\TR(a,b)$. Due to the critical exponents $I_\nu$ will not satisfy the Palais-Smale condition in general. Recall that solutions of \eqref{mainequ} satisfy the Pohozaev identity
\[
    P_\nu(u,v) := \int_\RN \sbr{|\nabla u|^2+|\nabla v|^2}
			  -\int_{\RN} \sbr{|u|^{2^*}+|v|^{2^*}} -\nu\ga_{\al+\beta}\int_{\RN}|u|^\al|v|^\beta=0.
\]
where 
\be\lab{gamma}
	\ga_p:=\f{N(p-2)}{2}.
\ee
Setting
\be\lab{pho}
	\PR_\nu(a,b):=\lbr{(u,v)\in \TR(a,b): P_\nu(u,v)=0},
\ee
we are interested in the problem whether
\be\lab{min}
	m_\nu(a,b):=\inf_{(u,v)\in\PR_\nu(a,b)} I_\nu(u,v).
\ee
is achieved. A {\it normalized ground state} is a solution $(u,v)\in\PR_\nu(a,b)$ of \eqref{mainequ} that achieves $m_\nu(a,b)$. We first see that a ground state can only exist if the interaction is attractive.

\bt\lab{thm2}
	If \eqref{assumption} holds and $\nu\le0$ then $m_\nu(a,b)$ is not achieved.
\et

Next we consider the case $\nu>0$ of attractive interaction. 

\bt\lab{thm1}
	Suppose \eqref{assumption} holds. 
	\begin{itemize}[fullwidth,itemindent=1em]
	\item[a)] If $\al+\beta<2+\f{4}{N}$ then there exists a normalized ground state, provided $0<\nu<\nu_0$ where $\nu_0$ is explicitly defined in \eqref{nu0} below. The normalized ground state is a local minimizer of $I_\nu$ on $\TR(a,b)$.
	\item[b)] If  $N=4$ and $\al+\beta=2+\f{4}{N}=3$ then there exists a normalized ground state, provided $0<\nu<\nu_0$ where $\nu_0$ is explicitly defined in  \eqref{nu0crit} below. The normalized ground state is of mountain pass type as critical point of $I_\nu$ on $\TR(a,b)$.
	\item[c)] If $\al+\beta>2+\f{4}{N}$, then there exists $\nu_1\ge 0$ such that there exists a normalized ground state if, and only if, $\nu>\nu_1$. Moreover, we have 
		\be\lab{tem6-7-1}
			\begin{cases}
			\nu_1=0 &\quad\text{if $N=4$, or $N=3$ and $|\al-\beta|>2$,}\\
			\nu_1>0 &\quad\text{if $N=3$, $\al\ge2$, $\beta\ge2$.}
			\end{cases}\
		\ee
		The normalized ground state is of mountain pass type as critical point of $I_\nu$ on $\TR(a,b)$.
	\item[d)] Both components of the normalized ground states in a), b) are positive, radially symmetric and radially decreasing. The corresponding Lagrange multipliers are positive:  $\la_1,\la_2>0$.
	\end{itemize}
\et

\br\lab{rem:thm1}
a) In the setting of Theorem~\ref{thm1}, it is an interesting open problem whether there exists a second solution of mountain pass type. The difficulty is to control the mountain pass value at a level where the Palais-Smale condition holds (in the space of radial $H^1$ functions).

b) Observe that in case $N=3$ the assumption $\al+\beta<2^*$ from \eqref{assumption} and $\al,\beta\ge2$ imply $|\al-\beta|<2$, hence the two cases in \eqref{tem6-7-1} exclude each other. However, there are still cases in Theorem \ref{thm1}~c), $N=3$, where we do not know whether $\nu_1=0$ or $\nu_1>0$. We conjecture that $\nu_1=0$ if $\min\{\al,\beta\}<2$.

c)	When $N=3$ and $\al+\beta=2+\f{4}{N}$, we can not obtain a similar result as Theorem \ref{thm1} b).
Indeed, in such case, we could obtain two numbers $\nu_0>0,\nu_1\ge0$ such that there exists a normalized ground state provided
both $\nu<\nu_0$ and $\nu>\nu_1$. However, we are unable to compare the values of $\nu_0$ and $\nu_1$, i.e., there may hold $\nu_0\le\nu_1$, while when $N=4$ and $\al+\beta=2+\f{4}{N}$, it acturally holds $\nu_1=0$ and $\nu_0>\nu_1$.
\er

An important part of the proof of Theorems~\ref{thm2} and \ref{thm1} will be to control $m_\nu(a,b)$. Since this has some interest in itself we collect our results in the next proposition. We need the best Sobolev embedding constant
\be\lab{Sobolev1}
	\SR:=\inf_{u\in\DR^{1,2}(\RN)\setminus\{0\}}\f{\nt{\nabla u}^2}{\ns{u}^2}.
\ee

\begin{proposition}\lab{prop1}
	Suppose \eqref{assumption} holds. 
	\begin{itemize}[fullwidth,itemindent=1em]
	\item[a)]	If $\al+\beta<2+\f{4}{N}$, then $m_\nu(a,b)$ is nonincreasing with respect to $\nu\in(-\iy,\nu_0)$ and
			$$m_\nu(a,b) \ \ 
				\begin{cases}
					=\f{1}{N}\SR^{N/2}   \quad &\text{when}~\nu\le0,\\
					<0					 &\text{when}~0<\nu<\nu_0.
				\end{cases}$$
	\item[b)]	If $\al+\beta=2+\f{4}{N}$ with $N=4$, then	 $m_\nu(a,b)$ is nonincreasing with respect to $\nu\in(-\iy,\nu_0)$, and
			$$m_\nu(a,b) \ \ 
				\begin{cases}
					=\f{1}{N}\SR^{N/2}		\quad &\text{when}~\nu\le0,\\
					<\f{1}{N}\SR^{N/2} 	\quad &\text{when}~0<\nu<\nu_0.
				\end{cases}$$
	\item[c)]	If $\al+\beta>2+\f{4}{N}$, then $m_\nu(a,b)$ is nonincreasing with respect to $\nu\in\R$, and
			$$m_\nu(a,b)	
				\begin{cases}
					=\f{1}{N}\SR^{N/2}		\quad &\text{when}~\nu\le\nu_1,\\
					<\f{1}{N}\SR^{N/2} 	\quad &\text{when}~\nu>\nu_1,\\
					\to0^+,					\quad &\text{when}~\nu\to+\iy.
				\end{cases}$$
	\end{itemize}
\end{proposition}
\br
a) In the case $\al+\beta<2+\f{4}{N}$ and $0<\nu<\nu_0$ the manifold $\PR_\nu(a,b)$ decomposes into two components $\PR_\nu^\pm(a,b)$; see section 2. Minimizing $I_\nu$ on $\PR_\nu^+(a,b)$ yields the normalized ground state from Theorem~\ref{thm1}~a). A minimizer of $I_\nu$ on the  $\PR_\nu^-(a,b)$, whose existence is an open problem, would be a critical point of $I_\nu$ of mountain pass type; cf. Remark~\ref{rem:thm1}~a). 

b) When $\al+\beta\le 2+\f{4}{N}$ the behavior of $m_\nu(a,b)$ in the range $\nu\ge\nu_0$ is an open problem; cf.\ also \cite{Soave=JFA=2020,Wei-Wu=2020}.
\er

Finally we  are concerned with asymptotic properties of the normalized ground states of \eqref{mainequ} obtained in Theorem \ref{thm1} as $\nu\to0^+$ or $\nu\to+\iy$. Similar results have been achieved for scalar equations in \cite{Wei-Wu=2020}. One type of limit behavior is related to the best Sobolev embedding constant $\SR$ defined in \eqref{Sobolev1}. It is well known that $\SR$ is achieved by the family
\be\lab{soliton}
  U_\e(x)=\big(N(N-2)\big)^{\f{N-2}{4}}\left(\f{\e}{\e^2+|x|^2}\right)^{\f{N-2}{2}} = \e^{\f{2-N}{2}}U_1(x/\e),\quad\e>0,
\ee
of radial functions, uniquely up to translations. For another type of limit behavior the following system appears as limit system:
\be\lab{limitequ}
	\left\{ 	\begin{aligned}
	&-\dl u+\la_1 u = {\al} |u|^{\al-2}|v|^\beta u &&\quad \text{in }\RN,\\
	&-\dl v+\la_2 v = {\beta} |u|^\al |v|^{\beta-2}v &&\quad \text{in }\RN,\\
	&\int_\RN u^2=a^2,\int_\RN v^2=b^2.
	\end{aligned}	\right.
\ee
Since $\al+\beta<2^*$ the problem \eqref{limitequ} is Sobolev subcritical. Normalized solutions of \eqref{limitequ} can be obtained as the critical points
of the energy functional $K(u,v): H^1(\RN)\times H^1(\RN) \to\R$
\be\lab{limitenergy}
	K(u,v):=\f{1}{2}\int_\RN \sbr{|\nabla u|^2+|\nabla v|^2} - \int_{\RN}|u|^\al|v|^\beta
\ee
on the constraint $\TR(a,b)$, with the parameters $\la_1,\la_2\in\R$ appearing as the Lagrange multipliers. The Pohozaev identity for solutions of \eqref{limitequ} is given by
\be
	L(u,v) := \int_\RN \sbr{|\nabla u|^2+|\nabla v|^2} -\ga_{\al+\beta}\int_{\RN}|u|^\al|v|^\beta = 0
\ee
with $\p$ defined by \eqref{gamma}. The investigation of \eqref{limitequ} has its own interest. A {\it normalized ground state} solution of \eqref{limitequ} is a solution $(u,v)$ that minimizes $K$ on
\be\lab{limitpho}
	\LR(a,b):=\lbr{(u,v)\in \TR(a,b): L(u,v)=0},
\ee
that is,
\be\lab{limitmin}
	K(u,v) =\inf_{(u,v)\in\LR(a,b)} K(u,v) =: l(a,b).
\ee

\bt\lab{thmlim}
	Suppose \eqref{assumption} holds and $\al+\beta \ne 2+\f{4}{N}$. Then $l(a,b)$ is achieved by a normalized ground state solution of \eqref{limitequ} with $\la_1,\la_2>0$. Both components of the solution are positive, radially symmetric and radially decreasing.
\et

By \cite[Theorem~1]{Busca-Sirakov:2000} every positive solution of \eqref{limitequ} with $\la_1,\la_2>0$ is radially symmetric with respect to some point.  Our next result relates the normalized ground states of \eqref{limitequ} to the vector-valued Gagliardo-Nirenberg inequality discussed in Appendix~A. Let $Z\in H^1(\R^N)$ be the unique positive radial solution of $-\dl w+w = w^{\al+\beta-1}$.

\begin{proposition}\lab{prop2}
	Suppose \eqref{assumption} holds.
	If $\al+\beta\neq2+\f{4}{N}$ and $\f{a^2}{b^2}=\f{\al}{\beta}$, then \eqref{limitequ}  with $\la_1,\la_2>0$ has a unique positive ground state solution $(W_1,W_2)$. This is of the form $W_i(x)=\sigma_iZ(\mu x)$, $i=1,2$, with $\sigma_1,\sigma_2,\mu$ explicitly given in \eqref{mu}-\eqref{sigma} below.
\end{proposition}

Motivated by \cite{Wei-Wu=2020} we also describe the asymptotic behavior of the normalized ground states of~\eqref{mainequ} from Theorem \ref{thm1}. We define
\be\lab{limgd}
\begin{aligned}
	\GR(a,b) := \big\{ (u,v)\in\TR(a,b): &\ (u,v)~\text{is a normalized ground state solution of~\eqref{limitequ}}\\
	                                                    &\ \text{with $\la_1,\la_2>0$}\big\}.
\end{aligned}
\ee
and 
\[
  \UR:=\lbr{U_\e:\e>0}
\]
with $U_\e$ from \eqref{soliton}, and we introduce the $L^2$-invariant scaling $s \s u (x) := e^{\frac N2 s}u(e^sx)$ and $s\s(u,v):=(s\s u,s\s v)$.

\bt\lab{thm3}
	Suppose \eqref{assumption} holds and let $(u_\nu,v_\nu)$ be a family of normalized ground state solutions of~\eqref{mainequ} from Theorem \ref{thm1} with $\nu\to0$. 
	\begin{itemize}[fullwidth,itemindent=2em]
	\item[a)]	If $\al+\beta<2+\f{4}{N}$ then there exists $t_\nu\in\R$ such that
			$$dist_H\Big(t_\nu\s(u_\nu,v_\nu),~\GR(a,b) \Big)\to0,\quad\text{as}~\nu\to0^+$$
		where $dist_H$ means the distance in $H=H^1(\RN)\times H^1(\RN)$. Moreover $t_\nu\sim (\p-2)^{-1}\ln \nu$.
	\item[b)]	If $\al+\beta\ge 2+\f{4}{N}$ then
			$$dist_{\DR^{1,2}}\Big( (u_\nu,v_\nu),~\UR\times\{0\} \Big)\to0,\quad\text{as}~\nu\to0^+,$$
		or
			$$dist_{\DR^{1,2}}\Big( (u_\nu,v_\nu),~\{0\}\times\UR \Big)\to0,\quad\text{as}~\nu\to0^+.$$
	\end{itemize}
\et

Recall from Theorem~\ref{thm1} that in the case $\al+\beta > 2+\f{4}{N}$ there need not exist a family $(u_\nu,v_\nu)$ with $\nu\to0$ as in Theorem~\ref{thm3}.

\bt\lab{thm4}
	Suppose \eqref{assumption} holds and $\al+\beta>2+\f{4}{N}$. Let $(u_\nu,v_\nu)$ be a family of normalized ground state solutions of~\eqref{mainequ} from Theorem \ref{thm1} with $\nu\to\infty$. Then
		$$dist_H\Big( s_\nu\s(u_\nu,v_\nu),~\GR(a,b) \Big)\to0,\quad\text{as}~\nu\to+\iy,$$
	where $s_\nu= (\p-2)^{-1}\ln \nu$.
\et

The paper is organized as follows. In Section 2 we prove our existence and non-existence results of normalized ground state solutions in the attractive case $\nu>0$. There we also investigate the behavior of the ground state energy as a function of $\nu>0$ and prove Proposition~\ref{prop1}. The non-existence result in Theorem~\ref{thm2} in the case $\nu\le0$ will be proved in Section 3. Section 4 contains the discussion of the limit system \eqref{limitequ}. Finally, in Section 5 we prove Theorems~\ref{thm3} and \ref{thm4} about the asymptotic behavior of the normalized ground state solutions as $\nu\to0$ or $\nu\to\infty$. The Appendix A contains some facts about the vector-valued Gagliardo-Nirenberg inequality. Condition \eqref{assumption} will be assumed throughout the paper.

\vskip0.3in
\section{Existence in the attractive case \texorpdfstring{$\nu>0$}{}}
	
	In this section we fix $\nu>0$, and start with the following compactness result.
\begin{proposition}\lab{PS2}
	Assume that
	\be\lab{monotonicity1}
		m_\nu(a,b)\le m_\nu(a_1,b_1)\quad\text{for any}~0<a_1\le a,~0<b_1\le b.
	\ee
	Let $\{(u_n,v_n)\}\subset \TR(a,b)$ be a sequence consisting of radially symmetric functions such that as $n\to+\iy$,
	\be\lab{tem3-1-1}
		I_\nu'(u_n,v_n)+\la_{1,n}u_n+\la_{2,n}v_n\to0\quad \text{for some}~\la_{1,n},\la_{2,n}\in\R,
	\ee
	\be\lab{tem3-2-2}
		I_\nu(u_n,v_n)\to c,\quad P_\nu(u_n,v_n)\to0,
	\ee
	\be\lab{tem3-3-3}
		u_n^-,v_n^-\to0,~\text{a.e. in}~\RN,
	\ee
with
\be\lab{levelcon}
	c\neq0\quad\text{and}\quad c<\f{1}{N}\SR^{N/2}+\min\lbr{ 0, m_\nu(a,b) }.
\ee
Then there exists a $(u,v)\in H_{rad}$, $u,v>0$ and $\la_1,\la_2>0$ such that
up to a subsequence $(u_n,v_n)\to(u,v)$ in $H$ and $(\la_{1,n},\la_{2,n})\to(\la_1,\la_2)$ in $\R^2$.
\end{proposition}

\bp The proof is divided into three steps.

\vskip0.1in
\noindent {\bf Step 1.} We show that $\{(u_n,v_n)\}$ is bounded in $H$, and $\la_{1,n},\la_{2,n}$ are bounded in $\R$. If $2+\f{4}{N}\le \al+\beta<2^*$, hence $ \p\ge2$, we obtain using $P_\nu(u_n,v_n)\to0$:
for $n$ large enough,
\begin{align*}
	c+1&\ge I_\nu(u_n,v_n)-\f{1}{2}P_\nu(u_n,v_n)\\
		&=\f{1}{N}(\ns{u_n}^{2^*}+\ns{v_n}^{2^*})+\f{\p-2}{2}\nu\int_\RN|u_n|^\al|v_n|^\beta \\
		&\ge C(\nt{\nabla u_n}^2+\nt{\nabla v_n}^2+o_n(1))
\end{align*}
for some $C>0$. This implies that $\{(u_n,v_n)\}$ is bounded in $H$. If $2<\al+\beta<2+\f{4}{N}$  we have $\p<2$, and obtain for $n$ large enough, using the inequality \eqref{GNine}:
\begin{align*}
	c+1&\ge I_\nu(u_n,v_n)-\f{1}{2^*}P_\nu(u_n,v_n)\\
		&=\f{1}{N}(\nt{\nabla u_n}^2+\nt{\nabla v_n}^2)-\f{2^*-\p}{2^*}\nu\int_\RN|u_n|^\al|v_n|^\beta \\
		&\ge \f{1}{N}(\nt{\nabla u_n}^2+\nt{\nabla v_n}^2)-C(\nt{\nabla u_n}^2+\nt{\nabla v_n}^2)^{\f{\p}{2}},
\end{align*}
for some $C>0$, which implies that $\{(u_n,v_n)\}$ is bounded in $H$. Moreover, from \eqref{tem3-1-1}, we observe that
	$$\la_{1,n}=-\f{1}{a^2}I_\nu'(u_n,v_n)[(u_n,0)]+o_n(1)	\quad\text{and}\quad
		\la_{2,n}=-\f{1}{b^2}I_\nu'(u_n,v_n)[(0,v_n)]+o_n(1), $$
so $\la_{1,n},\la_{2,n}$ are bounded also in $\R$.
Thus, there exist $(u,v)\in H_{rad}$, $\la_1,\la_2\in\R$ such that up to a subsequence,
	$$(u_n,v_n)\rh(u,v)\quad\text{in} ~H_{rad},$$
	$$(u_n,v_n)\ra(u,v)\quad\text{in} ~L^q(\RN)\times L^q(\RN),~\text{for}~2<q<2^*,$$
	$$(u_n,v_n)\ra(u,v)\quad\text{a.e. in}~\RN,$$
	$$(\la_{1,n},\la_{2,n})\to(\la_1,\la_2)\quad\text{in} ~\R^2.$$
Moreover, \eqref{tem3-1-1} and \eqref{tem3-3-3} imply
\be\lab{tem3-4-4}
	\left\{ 	\begin{aligned}
	&I_\nu'(u,v)+\la_1 u+\la_2v=0,\\
	&u\ge0,v\ge0,
	\end{aligned}\right.
\ee
and thus $P_\nu(u,v)=0$.

\vskip0.1in
\noindent {\bf Step 2.} We show that the weak limit satisfies $u\ne0$ and $v\ne0$, hence $u,v>0$ by the maximum principle. Arguing by contradiction we assume that $u=0$ so that \eqref{tem3-4-4} turns into
	$$-\dl v+\la_2 v=v^{2^*-1},\quad v\ge0,\quad v\in H^1(\RN).$$
If $v\not\equiv0$, then by the maximum principle $v>0$, hence $v=U_\e$ for some $\e>0$. This is a contradiction because $U_\e\not\in L^2(\RN)$, and therefore $v\equiv0$. Without loss of generality, we may assume that $\nt{\nabla u_n}^2\to h_1\ge0$ and $\nt{\nabla v_n}^2\to h_2\ge0$. We claim that $h_1+h_2>0$. Indeed, if $h_1=h_2=0$ then
	$$I_\nu(u_n,v_n)=\f{1}{N}(\nt{\nabla u_n}^2+\nt{\nabla v_n}^2)+\f{1}{2^*}P_\nu(u_n,v_n)+o_n(1)\to0,$$
as $n\to+\iy$, which is in contradiction with $c\neq0$. From $P_\nu(u_n,v_n)\to0$, one can see that
	$$h_1^2+h_2^2=\lim_{n\to\iy}\nt{\nabla u_n}^2+\nt{\nabla v_n}^2=\lim_{n\to\iy}\ns{u_n}^{2^*}+\ns{v_n}^{2^*}\le
		\SR^{-2^*/2}(h_1^{2^*}+h_2^{2^*}).$$
This implies
\[
\begin{aligned}
  c&=\lim_{n\to\iy}I_\nu(u_n,v_v)=\f{1}{N}(h_1^2+h_2^2)\\
   &\ge \f{1}{N}\min\lbr{ \eta_1^2+\eta_2^2:\,\eta_1^2+\eta_2^2\le \SR^{-2^*/2}(\eta_1^{2^*}+\eta_2^{2^*}),~\eta_1+\eta_2>0,~\eta_1\ge0, ~\eta_2\ge0}\\
		&=\f{1}{N}\min\lbr{ r^2:~r^2\le r^{2^*}\SR^{-2^*/2}\sbr{ \cos^{2^*}\q
			+\sin^{2^*}\q },~r\neq0,~\q\in[0,\f{\pi}{2}] }\\
		&=\f{1}{N}\SR^{N/2},
	\end{aligned}
\]
which is in contradiction with \eqref{levelcon}. So $u\neq0$, and analogously $v\neq0$.

\vskip0.1in
\noindent {\bf Step 3.} We show the strong convergence. Setting $(\bar u_n,\bar v_n)=(u_n-u,v_n-v)$ the Brezis-Lieb type lemma
\cite[Lemma 2.3]{Chen-Zou=TransAMS=2015} yields
	$$\int_\RN|u_n|^\al|v_n|^\beta-|u|^\al|v|^\beta-|\bar u_n|^\al|\bar v_n|^\beta=o_n(1).$$
Together with the weak convergence
	$$P_\nu(\bar u_n,\bar v_n)=P_\nu(u_n,v_n)-P_\nu(u,v)=o_n(1),$$
we deduce
	$$\lim_{n\to\iy}\big(\nt{\nabla\bar  u_n}^2+\nt{\nabla\bar  v_n}^2\big) = \lim_{n\to\iy}\big(\ns{\bar u_n}^{2^*}+\ns{\bar v_n}^{2^*}\big).$$
The same argument as in Step 2 gives that $\lim_{n\to\iy}\big(\nt{\nabla\bar  u_n}^2+\nt{\nabla\bar  v_n}^2\big) = 0$ or $\lim_{n\to\iy}\big(\nt{\nabla\bar  u_n}^2+\nt{\nabla\bar  v_n}^2\big) \ge \SR^{N/2}$. If the latter inequality holds we get
\[
	\begin{aligned}
	c&=\lim_{n\to\iy}I_\nu(u_n,v_v)=I_\nu(u,v)+\lim_{n\to\iy}I_\nu(\bar u_n,\bar v_n)\\
		&\ge m_\nu(\nt{u},\nt{v})+\lim_{n\to\iy}\f{1}{N}(\nt{\nabla\bar  u_n}^2+\nt{\nabla\bar  v_n}^2)\\
		&\ge m_\nu(a,b)+\f{1}{N}\SR^{N/2},
	\end{aligned}
\]
which contradicts \eqref{levelcon}. Here we used $0<\nt{u}\le a$, $0<\nt{v}\le b$ and assumption \eqref{monotonicity1}. Consequently we must have $\lim_{n\to\iy}\big(\nt{\nabla\bar  u_n}^2+\nt{\nabla\bar  v_n}^2\big) = 0$. We next claim that $\la_1,\la_2>0$. Indeed, if $\la_1\le0$, then
	$$-\dl u=|\la_1|u+u^{2^*-1}+\nu\al u^{\al-1}v^\beta\ge0,\quad\text{in}~\RN.$$
Then \cite[Lemma A.2]{Ikoma} implies $u=0$, which is a contradiction. Therefore $\la_1>0$, and analogously $\la_2>0$, as claimed. Combining \eqref{tem3-1-1}, \eqref{tem3-2-2}, \eqref{tem3-4-4} and $P_\nu(u,v)=0$, we obtain
\be\lab{tem3-5}
	\la_1 a^2+\la_2 b^2=\la_1 \nt{u}^2+\la_2\nt{v}^2.
\ee
It follows that $\nt{u}=a$, $\nt{v}=b$ and hence $(u_n,v_n)\to(u,v)$ in $H$.
\ep

\vskip0.1in
\subsection{The case \texorpdfstring{$\al+\beta<2+\f{4}{N}$}{}}

We define
\be\lab{nu0}
	\nu_0 := \min\lbr{\f{1}{\p},\f{2^*+2-\p}{2\cdot2^*}}\cdot\f{ \SR^{\f{2^*(2-\p)}{2(2^*-2)}} }{\CR(N,\al,\beta)}
		 \cdot\f{ (2^*-2)(2-\p)^{\f{2-\p}{2^*-2}} }{ (2^*-\p)^{\f{2^*-\p}{2^*-2}} (a^2+b^2)^{\f{(\al+\beta-\p)}{2}}}
\ee
where $\CR(N,\al,\beta)$ is defined in \eqref{VGNconstant}, which comes from the vector-valued version of the Gagliardo-Nirenberg inequality. For $\nu>0$ and $(u,v)\in\TR(a,b)$ we consider the map $\Phi^\nu_{(u,v)}:\R\to\R$ defined by
	$$\Phi^\nu_{(u,v)}(t) := I_\nu\big(t\s(u,v)\big) 
	      = \f{1}{2}e^{2t}\int_\RN \sbr{|\nabla u|^2+|\nabla v|^2}-\f{1}{2^*}e^{2^*t}\int_{\RN} \sbr{|u|^p+|v|^p}-\nu e^{\p t}
					\int_{\RN}|u|^\al|v|^\beta.$$
Observe that $(\Phi^\nu_{(u,v)})'(t)=P_\nu(t\s(u,v))$, hence
	$$\PR_\nu(a,b):=\lbr{(u,v)\in \TR(a,b): (\Phi^\nu_{(u,v)})'(0)=0 }.$$
We decompose $\PR_\nu(a,b)$ into three disjoint sets
 $\PR_\nu(a,b)=\PR_\nu^+(a,b)\cup\PR_\nu^0(a,b)\cup\PR_\nu^-(a,b)$,
defined by
\begin{align*}
	\PR_\nu^+(a,b) &:=\left\{ (u,v)\in \PR_\nu(a,b):(\Phi^\nu_{(u,v)})''(0)>0 \right\},\\
	\PR_\nu^0(a,b) &:=\left\{ (u,v)\in \PR_\nu(a,b):(\Phi^\nu_{(u,v)})''(0)=0 \right\},\\
	\PR_\nu^-(a,b) &:=\left\{ (u,v)\in \PR_\nu(a,b):(\Phi^\nu_{(u,v)})''(0)<0 \right\}.
\end{align*}
For $\rho>0$ we define
\[
	h(\rho):=\f{1}{2}\rho^2-A\rho^{\p}-B\rho^{2^*}\quad\text{and}\quad g(\rho):=\rho^{2-\p}-2^*B\rho^{2^*-\p}
\]
with
	$$A:=\nu\CR(N,\al,\beta)(a^2+b^2)^{\f{\al+\beta-\p}{2}}\quad\text{and}\quad B:=\f{1}{2^*}\SR^{-2^*/2}.$$
Then using \eqref{VGNine} we obtain for any $u,v\in H^1(\RN)$:
	$$I_\nu(u,v)\ge h\sbr{\big(\nt{\nabla u}^2+\nt{\nabla v}^2\big)^{1/2}}$$
and
	$$h'(p)=\rho^{\p-1}\sbr{g(\rho)-\p A}.$$
Setting
\be\lab{rho}
	\rho_*:=\sbr{\f{2-\p}{2^*(2^*-\p)B}}^{1/(2^*-2)}
\ee
it is easy to check that $g(\rho)$ is strictly increasing in $(0,\rho_*)$, and is strictly decreasing in $(\rho_*,+\iy)$. By direct computations, 
the assumption $0<\nu<\nu_0$ with $\nu_0$ defined in \eqref{nu0} implies $g(\rho_*)>\p A$ and $h(\rho_*)>0$. This means that 
$h(\rho)$ has exactly two critical points $0<\rho_1<\rho_*<\rho_2$ with
	$$h(\rho_1)=\min_{0<\rho<\rho_*}h(\rho)<0\quad\text{and}\quad h(\rho_2)=\max_{\rho>0} h(\rho)>0.$$
Moreover, there exist $R_1>R_0>0$ such that $h(R_0)=h(R_1)=0$ and $h(\rho)>0$ iff $\rho\in(R_0,R_1)$.

\bl\lab{structure1}
	Let $0<\nu<\nu_0$. Then $\Phi^\nu_{(u,v)}$ has, for every $(u,v)\in \TR(a,b)$, exactly two critical	points $t_\nu(u,v)<s_\nu(u,v)$ 
	and two zeros $c_\nu(u,v)<d_\nu(u,v)$ satisfying $t_\nu(u,v)<c_\nu(u,v)<s_\nu(u,v)<d_\nu(u,v)$. Moreover:
	\begin{itemize}[fullwidth,itemindent=1em]
		\item[a)]	$\PR_\nu^0(a,b)=\emptyset$ and $\PR_\nu(a,b)$ is a submanifold of $H$.
		\item[b)]	$t\s(u,v)\in\PR_\nu^+(a,\nu)$ if and only if $t=t_\nu(u,v)$, and $t\s(u,v)\in\PR_\nu^-(a,b)$ if and only if $t=s_\nu(u,v)$.
		\item[c)] 	$\sbr{\nt{\nabla (t\s u)}^2+\nt{\nabla (t\s v)}^2}^{1/2}\le R_0$ for every $t\le c_\nu(u,v)$ and
				\be
	\begin{aligned}
					I_\nu(t_\nu(u,v)\s(u,v))&=\min\lbr{ I_\nu(t\s(u,v)): t\in\R ~\text{and}~
						\sbr{\nt{\nabla (t\s u)}^2+\nt{\nabla (t\s v)}^2}^{1/2}\le R_0}\\
&<0.
\end{aligned}
				\ee
		\item[d)] 	$\Phi^\nu_{(u,v)}(t)$ is strictly decreasing on $(s_\nu(u,v),+\iy)$, and
					$\Phi^\nu_{(u,v)}(s_\nu\s(u,v))=\max\limits_{t\in\R}\Phi^\nu_{(u,v)}(t)>0$.
		\item[e)] 	The maps $(u,v)\mapsto t_\nu(u,v)$ and $(u,v)\mapsto s_\nu(u,v)$ are of class $C^1$.
	\end{itemize}
\el

\bp
	We prove it briefly, a similar proof can be found in \cite{Li-Zou-1}. For any $(u,v)\in \TR(a,b)$ we have
	$$\Phi_{(u,v)}^\nu(s)=I_\nu(s\s(u,v))\geq h\sbr{e^s(\nt{\nabla u}^2+\nt{\nabla v}^2)^{1/2}}$$
and
	$$\Phi_{(u,v)}^\nu(s)>0,\quad \text{for all }~s\in \sbr{ \log\f{R_0}{(\nt{\nabla u}^2+\nt{\nabla v}^2)^{1/2}},~
		\log\f{R_1}{(\nt{\nabla u}^2+\nt{\nabla v}^2)^{1/2} } }$$
because $0<\nu<\nu_0$. Recalling the facts that $\Phi_{(u,v)}^\nu(-\infty)=0^-$ and $\Phi_{(u,v)}^\nu(+\infty)=-\infty$, 
we see that $\Phi_{(u,v)}^\nu(s)$ has at least two critical points $t_\nu(u,v)<s_\nu(u,v)$, where $s_\nu(u,v)$ is the global 
maximum point at positive level, and $t_\nu(u,v)$ is the minimum point on $\sbr{-\iy,\log\f{R_0}{(\nt{\nabla u}^2+\nt{\nabla v}^2)^{1/2}}}$ 
at negative level. On the other hand, $\Phi_{(u,v)}^\nu(s)$ has at most two critical points in $\R$ by \cite[Remark 3.1]{Li-Zou-1}, 
hence it has exactly the two critical points $t_\nu(u,v)$ and $s_\nu(u,v)$. Now $s\s(u,v)\in\PR_\nu(a,b)$ implies $s=t_\nu(u,v)$ 
or $s=s_\nu(u,v)$ because $(\Phi_{(u,v)}^\nu)'(s)=P_\nu(s\s(u,v))$. Moreover, as in \cite[Lemma 3.2]{Li-Zou-1}, $0<\nu<\nu_0$
 implies $\PR_\nu^0(a,b)=\emptyset$ and $\PR_\nu(a,b)$ is a submanifold of $H$. It follows that $t_\nu(u,v)\s(u,v)\in\PR_\nu^+(a,b)$ 
 and $s_\nu(u,v)\s(u,v)\in\PR_\nu^-(a,b)$. By the monotonicity, $\Phi_{(u,v)}^\nu(s)$ has exactly two zeros $c_\nu(u,v)$ and $d_\nu(u,v)$, 
 with $t_\nu(u,v)<c_\nu(u,v)<s_\nu(u,v)<d_\nu(u,v)$. Finally, the implicit function theorem implies that the maps $(u,v)\mapsto t_\nu(u,v)$
  and $(u,v)\mapsto s_\nu(u,v)$ are of class $C^1$.
\ep

For $R>0$ we define
	$$A_R(a,b):=\lbr{(u,v)\in \TR(a,b):\sbr{\nt{\nabla u}^2+\nt{\nabla v}^2}^{1/2}<R }.$$
	
\bl
	If $0<\nu<\nu_0$ then the following statements hold.
	\begin{itemize}[fullwidth,itemindent=1em]
	\item[a)] 	$m_\nu(a,b)=\inf\limits_{A_{R_0}(a,b)}I_\nu(u,v)<0$
	\item[b)]	$m_\nu(a,b)\le m_\nu(a_1,b_1)$ for any $0<a_1\le a$, $0<b_1\le b$
	\end{itemize}	
\el

\bp
a)
Lemma \ref{structure1} implies $\PR_\nu^+(a,b)\subset A_{R_0}(a,b)$ and
	$$m_\nu(a,b)=\inf_{\PR_\nu(a,b)}I(u,v)=\inf_{\PR_\nu^+(a,b)} I(u,v)<0.$$
Obviously, $m_\nu(a,b)\ge\inf_{A_{R_0}(a,b)}I_\nu(u,v)$. On the other hand, for any $(u,v)\in A_{R_0}(a,b)$, we have
	$$m_\nu(a,b)\leq I_\nu(t_\nu(u,v)\s(u,v))\le I_\nu(u,v).$$
It follows that $m_\nu(a,b)=\inf\limits_{A_{R_0}(a,b)}I_\nu(u,v)$.

b)
This can be proved with minor modifications by following the strategy in \cite[Lemma 3.2]{Jeanjean-Lu=CVPDE=2020}, where a scalar 
equation is considered. Recalling the definition of $\rho_*$ in \eqref{rho} we obtain as in a)
\be\lab{tem3-6-7}
	m_\nu(a,b)=\inf_{A_{\rho_*}(a,b)} I_\nu(u,v).
\ee
We prove $m_\nu(a,b)\le m_\nu(a_1,b_1)+\e$ for arbitrary $\e>0$. Let $(u,v)\in A_{\rho_*}(a_1,b_1)$ be such that 
$I_\nu(u,v)\le m_\nu(a_1,b_1)+\f{\e}{2}$ and $\phi\in\CR_0^\iy(\RN)$ be a cut-off function satisfying
	$$0\le\phi\le1\quad\text{and}\quad \phi(x)=\begin{cases}0\quad\text{if}~|x|\ge2,\\1\quad\text{if}~|x|\le1. \end{cases}$$
For $\delta>0$ we consider $u_\delta(x)=u(x)\phi(\delta x)$ and $v_\delta(x)=v(x)\phi(\delta x)$; clearly $(u_\delta,v_\delta)\to(u,v)$ 
in $H$ as $\delta\to0$. As a consequence, for $\eta>0$ small there exists $\delta>0$ small such that
\be\lab{tem31}
	I_\nu(u_\delta,v_\delta)\le I_\nu(u,v)+\f{\e}{4}\quad\text{and}\quad
			\sbr{\nt{\nabla u_\delta}^2+\nt{\nabla v_\delta}^2}^{1/2}<\rho_*-\eta.
\ee
Now we take $\vp\in\CR_0^\iy(\RN)$ such that $\supp(\vp)\subset B(0,1+\f{4}{\delta})\setminus B(0,\f{4}{\delta})$,
where $B(0,r)$ is the ball with a radius   $r$   centered  at the origin.
Set
	$$w_a=\f{\sqrt{a^2-\nt{u_\delta}^2}}{\nt{\vp}}\cdot\vp \quad\text{and}\quad w_b=\f{\sqrt{b^2-\nt{u_\delta}^2}}{\nt{\vp}}\cdot\vp.$$
and observe that
	$$\sbr{supp(u_\delta)\cup supp(v_\delta)}\cap \sbr{supp(s\s w_a)\cup supp(s\s w_b)}=\emptyset,$$
for $s<0$, hence $(u_\delta+s\s w_a,v_\delta+s\s w_b)\in\TR(a,b)$. Next we have
\be\lab{tem32}
	I_\nu(s\s (w_a,w_b))\le\f{\e}{4}\quad\text{and}\quad
		\sbr{\nt{\nabla s\s w_a}^2+\nt{\nabla s\s w_b}^2}^{1/2}\le \f{\eta}{2}
\ee
for $s<0$ sufficiently close to $-\iy$ because $I_\nu\big(s\s (w_a,w_b)\big)\to0$ and
 $\sbr{\nt{\nabla s\s w_a}^2+\nt{\nabla s\s w_b}^2}^{1/2}\to0$ as $s\to-\iy$. It follows that
	$$\sbr{\nt{\nabla (u_\delta+s\s w_a)}^2+\nt{\nabla (v_\delta+s\s w_b)}^2}^{1/2}< \rho_*.$$
Now using \eqref{tem31} and \eqref{tem32}, we obtain
	$$m_\nu(a,b)\le I_\nu(u_\delta+s\s w_a,v_\delta+s\s w_b)
		=I_\nu(u_\delta,v_\delta)+I_\nu(s\s w_a,s\s w_b)\le m_\nu(a_1,b_1)+\e$$
which finishes the proof.
\ep

\vskip0.12in
\bp[Proof of Theorem \ref{thm1} a) and d)]
Let $(\tilde u_n,\tilde v_n)\in A_{R_0}(a,b)$ be a a minimizing sequence $m_\nu(a,b)$. By the symmetric decreasing rearrangement
we may assume that  $\tilde u_n$, $\tilde v_n$ are both radial. After passing to $(|\ti u_n|,|\ti v_n|)$ we may also assume that 
$(\tilde u_n,\tilde v_n)$ are nonnegative. Furthermore, using $I_\nu(t_+(\ti u_n,\ti v_n)\s(\ti u_n,\ti v_n))\le I_\nu(\ti u_n,\ti v_n)$ 
we can assume that $(\ti u_n,\ti v_n)\in\PR_\nu^+(a,b)$ for $n\ge1$. Therefore, by Ekeland's varational principle there is a 
radially symmetric Palais-Smale sequence $(u_n,v_n)$ for $I_\nu|_{\TR(a,b)}$ satisfying 
$||(u_n,v_n)-(\tilde u_n,\tilde v_n)||_H\to0$ as $n\to\infty$. This implies
	$$P_\nu(u_n,v_n)=P_\nu(\ti u_n,\ti v_n)+o(1)\to0\quad \text{and}\quad u_n^-,v_n^-\ra0\ \text{a.e. in}\ \RN.$$
Now Proposition \ref{PS2} with $c=m_\nu(a,b)$ implies that $(u_n,v_n)\to(u,v)$ in $H$ and $(\la_{1,n},\la_{2,n})\to(\la_1,\la_2)\in(\R^+)^2$ 
along a subsequence. By the strong convergence, $(u,v) \in \PR_\nu(a,b)$ is a solution of \eqref{mainequ} and hence a normalized ground state.
\ep

\vskip0.1in
\subsection{The case \texorpdfstring{$\al+\beta> 2+\f{4}{N}$}{}}\lab{tem6-11}

\bl\lab{structure2}
	Fix $\nu>0$. Then $\Phi^\nu_{(u,v)}$ has exactly one critical point $t_\nu(u,v)$ for every $(u,v)\in \TR(a,b)$. Moreover:
	\begin{itemize}[fullwidth,itemindent=1em]
		\item[a)]	$\PR_\nu(a,b)=\PR_\nu^-(a,b)$ and $\PR_\nu(a,b)$ is a submanifold of $H$.
		\item[b)]	$t\s(u,v)\in\PR_\nu(a,b)$ if and only if $t=t_\nu(u,v)$.
		\item[c)] 	$\Phi^\nu_{(u,v)}(t)$ is strictly dereasing and concave on $(t_\nu(u,v),+\iy)$ and
						$$\Phi^\nu_{(u,v)}(t_\nu(u,v))=\max_{t\in\R}\Phi^\nu_{(u,v)}(t)>0.$$
		\item[d)] 	The map $(u,v)\mapsto t_\nu(u,v)$ is of class $C^1$.
	\end{itemize}
\el
\bp
	For any $(u,v)\in\PR_\nu(a,b)\setminus\PR_\nu^-(a,b)$ there holds
		$$\nt{\nabla u}^2+\nt{\nabla v}^2=\ns{u}^{2^*}+\ns{v}^{2^*}+\p \nu\int_\RN|u|^\la|v|^\beta,$$
and
		$$2(\nt{\nabla u}^2+\nt{\nabla v}^2)\ge2^*(\ns{u}^{2^*}+\ns{v}^{2^*})+\p^2 \nu\int_\RN|u|^\la|v|^\beta.$$
This implies
	$$0\ge (2^*-2)(\ns{u}^{2^*}+\ns{v}^{2^*})+(\p-2)\p\nu\int_\RN|u|^\la|v|^\beta$$
which is impossible because $\p>2$. Consequently $\PR_\nu(a,b)=\PR_\nu^-(a,b)$ and as in \cite[Lemma 3.2]{Li-Zou-1} one sees that 
$\PR_\nu(a,b)$ is a submanifold of $H$. The conclusions $b)$-$d)$ follow from the properties of the fiber map $\Phi_{(u,v)}^\nu(s)$,
which can be proved as \cite[Lemma 6.2]{Soave=JDE=2020}. 
\ep

\bl\lab{minimax}
	The minimum $m_\nu(a,b)$ has a minimax description:
		$$m_\nu(a,b)=\inf_{\TR(a,b)}\max_{t\in\R}I_\nu(t\s(u,v)).$$
\el

\bp
	For $(u,v)\in\TR(a,b)$ Lemma \ref{structure2} implies
	$$\max_{t\in\R}I_\nu\big(t\s(u,v)\big) = I_\nu\big(t_\nu(u,v)\s(u,v)\big) \ge m_\nu(a,b),$$
hence $m_\nu(a,b)\le\inf\limits_{\TR(a,b)}\max\limits_{t\in\R}I_\nu\big(t\s(u,v)\big)$.
On the other hand, Lemma \ref{structure2} also implies
	$$I_\nu(u,v)=\max_{t\in\R}I_\nu\big(t\s(u,v)\big) \ge \inf_{\TR(a,b)}\max_{t\in\R}I_\nu\big(t\s(u,v)\big)\quad\text{if $(u,v)\in\PR_\nu(a,b)$,}$$
and thus $m_\nu(a,b) = \inf\limits_{\TR(a,b)}\max\limits_{t\in\R}I_\nu(t\s(u,v))$.
\ep

\bl\lab{PSseq}
	There exists a Palais-Smale sequence $(u_n,v_n)\subset \TR(a,b)\cap H_{rad}$ for $I_\nu|_{\TR(a,b)}$ at the
	level $m_\nu(a,b)$, with $P_\nu(u_n,v_n)\to0$ and $u_n^-,v_n^-\to0$ a.e. in $\RN$ as $n\to\iy$.
\el

\bp
	Setting $\tilde I_\nu(s,(u,v)):=I_\nu(s\s (u,v))$ Lemma~\ref{minimax} we see that 
	\[
	  m_\nu(a,b)=\inf_{\ga\in\Ga_\nu}\max_{t\in[0,1]} \tilde I_\nu(\ga(t))
	\]
	 where
	 $$\Ga_\nu = \lbr{\ga:[0,1]\to\R\times \big(\TR(a,b)\cap H_{rad}\big):
	       \text{$\ga$ is continuous, $\ga(0)\in\{0\}\times A_r(a,b)$, $\ga(1)\in\{0\}\times I_\nu^0$}}$$
for $r>0$ arbitrarily small. Then a required Palais-Smale sequence can be found in a standard way as in \cite{Bartsch-Jeanjean-Soave=JMPA=2016}.
\ep

	In view of Proposition \ref{PS2} we need some properties about $m_\nu(a,b)$ in order to obtain the compactness of the Palais-Smale sequence.

\bl\lab{property1}
	For $a,b>0$ fixed the following statements hold:
	\begin{itemize}
	\item[a)]	$m_\nu(a,b)\le m_\nu(a_1,b_1)$ for any $0< a_1\le a$ and $0<b_1\le b$.
	\item[b)]	$m_\nu(a,b)$ is nonincreasing with respect to $\nu\in(0,+\iy)$.
	\item[c)]	$\lim_{\nu\to+\iy} m_\nu(a,b)=0^+$
	\end{itemize}
\el

\bp
a) This can be proved as Lemma \ref{lem2-2}~b).\\

b) For any $\nu\ge\nu'>0$, we see that
		$$m_\nu(a,b) = \inf_{\TR(a,b)}\max_{t\in\R}I_\nu\big(t\s(u,v)\big)
			\le \inf_{\TR(a,b)}\max_{t\in\R}I_{\nu'}\big(t\s(u,v)\big) = m_{\nu'}(a,b).$$
Consequently $m_\nu(a,b)$ is nonincreasing with respect to $\nu>0$.\\

c) We first prove $m_\nu(a,b)>0$ for any $\nu>0$. Indeed, for $(u,v)\in\PR_\nu(a,b)$ there holds
	$$\begin{aligned}
	\nt{\nabla u}^2+\nt{\nabla v}^2&=\ns{u}^{2^*}+\ns{v}^{2^*}+\nu\p\int_\RN|u|^\al|v|^\beta\\
		&\le C_1\big(\nt{\nabla u}^2+\nt{\nabla v}^2\big)^{\f{2^*}{2}}+C_2\big(\nt{\nabla u}^2+\nt{\nabla v}^2\big)^{\f{\p}{2}}
	\end{aligned}$$
which implies $\inf\limits_{\PR_\nu(a,b)}\big(\nt{\nabla u}^2+\nt{\nabla v}^2\big) >0$. So we have
	$$\begin{aligned}
	m_\nu(a,b)&=\inf_{\PR_\nu(a,b)}I_\nu(u,v)=\inf_{\PR_\nu(a,b)} I_\nu(u,v)-\f{1}{\p}P_\nu(u,v)\\
	&=\inf_{\PR_\nu(a,b)} \f{\p-2}{2\p}\big(\nt{\nabla u}^2+\nt{\nabla v}^2\big)+\f{2^*-\p}{2^*\p}\nu\int_\RN|u|^\al|v|^\beta\\
	&\ge C\inf_{\PR_\nu(a,b)}\big(\nt{\nabla u}^2+\nt{\nabla v}^2\big) >0.
	\end{aligned}$$
	Now c) follows if we can show that for any $\e>0$ there exists $\bar\nu>0$ such that
\be\lab{tem6-7-4}
	m_\nu(a,b)<\e\quad\text{for any}~\nu\ge\bar\nu.
\ee
Choose $\phi\in\CR_0^\iy(\RN)$ with $\nt{\phi}\le \min\lbr{a,b}$. Then
\be\lab{tem6-7-3}
	m_\nu(a,b)\le m_\nu(\nt{\phi},\nt{\phi})\le \max_{t\in\R}I_\nu(t\s\phi,t\s\phi)
			=\max_{t\in\R} \sbr{2E(t\s\phi)-\nu e^{\p t}\nm{\phi}_{\al+\beta}^{\al+\beta}},
\ee
where $E(u)$ is defined by
\be\lab{defofE}
	E(u)=\f{1}{2}\int|\nabla u|^2-\f{1}{2^*}\int |u|^{2^*}.
\ee
It is easy to see that $E(t\s\phi)\to0^+$ as $t\to-\iy$, so there exists $t_0>0$ such that
$E(t\s\phi)<\f{\e}{4}$ for any $t<-t_0$. On the other hand, there exists  $\bar\nu > 0$ such that
	$$\max_{t\ge-t_0} \sbr{2E(t\s\phi)-\nu e^{\p t}\nm{\phi}_{\al+\beta}^{\al+\beta}}
		\le \f{2\nt{\nabla\phi}^N}{N\ns{\phi}^N}- \nu e^{-\p t_0}\nm{\phi}_{\al+\beta}^{\al+\beta} <\e\quad\text{for $\nu\ge \bar\nu$.}$$
It follows that $\displaystyle \max_{t\in\R} \sbr{2E(t\s\phi)-\nu e^{\p t}\nm{\phi}_{\al+\beta}^{\al+\beta}}<\e$ when $\nu\ge \bar\nu$.
Then combining with \eqref{tem6-7-3}, we obtain \eqref{tem6-7-4}.
\ep


	Lemma \ref{property1}~c) implies 
\be
	\nu_1:=\inf\lbr{ \nu>0:~m_\nu(a,b)<\f{1}{N}\SR^{N/2} } < +\iy.
\ee
In order to investigate in which case $\nu_1=0$ or $\nu_1>0$ we use the test functions $\phi U_\e$
where $U_\e$ is defined by \eqref{soliton} and $\phi\in\CR_0^\iy(\RN)$ is a radial cut-off function satisfying
	$$0\le\phi\le1\quad\text{and}\quad \phi(x)=\begin{cases}1\quad\text{if}~|x|\le1,\\0\quad\text{if}~|x|\ge2. \end{cases}$$
Then the following estimates hold; proofs can be found in \cite{Jeanjean=2020,Soave=JFA=2020}.

\bl\lab{est}
	For $\eta_\e=\phi U_\e$ we have as $\e\to0^+$:
		$$\nt{\nabla \eta_\e}^2=\SR^{N/2}+O(\e^{N-2}),\quad \ns{\eta_\e}^{2^*}=\SR^{N/2}+O(\e^N),$$
		$$\nt{\eta_\e}^2=\begin{cases} 	O(\e^2) 	\quad &\text{if}\quad N\ge5,\\
			O(\e^2|\log \e|) 						\quad&\text{if}\quad N=4,\\
			O(\e)									\quad&\text{if}\quad N=3,		
		\end{cases}$$
		$$\nm{\eta_\e}_p^p=\begin{cases}  O(\e^{N-(N-2)p/2})	\quad&\text{if}\quad N\ge4~\text{and}~p\in(2,2^*)~
								\text{or if}~N=3~\text{and}~p\in(3,6),\\
					O(\e^{p/2})					\quad&\text{if}\quad N=3~\text{and}~p\in(2,3),\\
					O(\e^{3/2}|\log\e|)			\quad&\text{if}\quad N=3~\text{and}~p=3.		\end{cases} $$
\el

\bl\lab{property2}
	For fixed $a,b>0$ the following statements hold.
	\begin{itemize}
	\item[a)]	$m_\nu(a,b)\le\f{1}{N}\SR^{N/2}$ for all $\nu>0$.
	\item[b)]	If $N=3$ and $|\al-\beta|>2$, or $N=4$, then $\nu_1=0$.
	\item[c)]	If $N=3$ and $\al,\beta\ge2$ then $\nu_1>0$.
	\item[d)]	If $0<\nu\le\nu_1$ then $m_\nu(a,b)=\f{1}{N}\SR^{N/2}$ is not achieved.
	\end{itemize}
\el

\bp
a) We may assume $\al>\beta$. Then $\al-\beta>2$ implies $\beta<2$ in dimension $N=3$ because $\al+\beta<6$. In dimension $N=4$ we have $\beta<2$ anyway because $\al+\beta<4$. Now we choose $\q>0$ such that
	$$\begin{cases} \q\in\big(\f{\al+\beta-2}{4\beta}\,,\,\f{6-\al-\beta}{4(2-\beta)}\big),\quad&\text{when}~N=3,\\
			\q=1,\quad&\text{when}~N=4.\end{cases}$$
We set $u_\e=\f{ a}{\nt{\eta_\e}}\eta_\e$ and $v_\e=\f{\e^\q}{\nt{\eta_\e}}\eta_\e$, and let $t_\e:=t_\nu(u_\e,v_\e)$ be given by Lemma \ref{structure2}. Then $t_\e\s(u_\e,v_\e)\in\PR_\nu(a,\e^\q)$, hence for $\e>0$ small we obtain:
\be\lab{tem4-13-8}
	\begin{aligned}
	m_\nu(a,b)&\le m_\nu(a,\e^\q)\le I_\nu\big(t_\e\s(u_\e,v_\e)\big)\\
		&=\f{1}{2}e^{2t_\e}\sbr{\nt{\nabla u_\e}^2+\nt{\nabla v_\e}^2}
				-\f{1}{2^*}e^{2^* t_\e}\sbr{\ns{u_\e}^{2^*}+\ns{v_\e}^{2^*}}
				-\nu e^{\p t_\e}\int_\RN u_\e^\al v_\e^\beta\\
		&\le \max_{s>0}\sbr{ \f{a^2+\e^{2\q}}{2}s^2\f{\nt{\nabla\eta_\e}^2}{\nt{\eta_\e}^2}
			-\f{a^{2^*}+\e^{2^*\q}}{2^*}s^{2^*}\f{\ns{\eta_\e}^{2^*}}{\nt{\eta_\e}^{2^*}} }
			-C\e^{\beta\q} e^{\p t_\e} \f{\nm{\eta_\e}_{\al+\beta}^{\al+\beta}}{\nt{\eta_\e}^{\al+\beta}}\\
		&=\f{1}{N}\sbr{ \f{(a^2+\e^{2\q})\nt{\nabla\eta_\e}^2 }
			{\sbr{ a^{2^*}+\e^{2^*\q} }^{2/2^*}\ns{\eta_\e}^2 }}^{N/2}
			-C\e^{\beta\q} e^{\p t_\e} \f{\nm{\eta_\e}_{\al+\beta}^{\al+\beta}}{\nt{\eta_\e}^{\al+\beta}}\\
		&=\f{1}{N}\SR^{N/2}+O(\e^{N-2})+O(\e^{2\q})
			-C\e^{\beta\q} e^{\p t_\e} \f{\nm{\eta_\e}_{\al+\beta}^{\al+\beta}}{\nt{\eta_\e}^{\al+\beta}}.
	\end{aligned}
\ee
	We claim that $e^{t_\e} \ge C\nt{\eta_\e}$ as $\e\to0$ for some constant $C>0$. By  definition $P_\nu\big(t_\e\s(u_\e,v_\e)\big)=0$ and thus
	$$e^{(2^*-2)t_\e}(\ns{u_\e}^{2^*}+\ns{v_\e}^{2^*})\le\nt{\nabla u_\e}^2+\nt{\nabla v_\e}^2.$$
It follows that
\be\lab{tem4-13-9}
	e^{t_\e}\le\sbr{\f{\nt{\nabla u_\e}^2+\nt{\nabla v_\e}^2}{\ns{u_\e}^{2^*}+\ns{v_\e}^{2^*}}}^{1/(2^*-2)}.
\ee
From \eqref{tem4-13-9} and $\p>2$, we deduce for $\e>0$ small:
\be\lab{tem4-13-7}
	\begin{aligned}
	&e^{(2^*-2)t_\e}
                 =\f{\nt{\nabla u_\e}^2+\nt{\nabla v_\e}^2}{\ns{u_\e}^{2^*}+\ns{v_\e}^{2^*}}
			-\p\f{ \nu\int_\RN u_\e^\al v_\e^\beta }{ \ns{u_\e}^{2^*}+\ns{v_\e}^{2^*} }e^{(\p-2)t_\e}\\
	&\hspace{.5cm}\ge \f{\nt{\nabla u_\e}^2+\nt{\nabla v_\e}^2}{\ns{u_\e}^{2^*}+\ns{v_\e}^{2^*}}
			-\p\f{\nu\int_\RN u_\e^\al v_\e^\beta}{ \ns{u_\e}^{2^*}+\ns{v_\e}^{2^*} }
			\sbr{\f{\nt{\nabla u_\e}^2+\nt{\nabla v_\e}^2}{\ns{u_\e}^{2^*}+\ns{v_\e}^{2^*}}}^{\f{\p-2}{2^*-2}}\\
	&\hspace{.5cm}= \f{\nt{\nabla\eta_\e}^2\nt{\eta_\e}^{2^*-2}}{\ns{\eta_\e}^{2^*}}
			\sbr{C_1-C_2\e^{\beta\q}\nt{\nabla\eta_\e}^{-2\f{2^*-\p}{2^*-2}}
			\ns{\eta_\e}^{2^*\f{\p-2}{2^*-2}}\f{\nm{\eta_\e}_{\al+\beta}^{\al+\beta}}{\nt{\eta_\e}^{\al+\beta-\p}}}
	\end{aligned}
\ee
Using $\nt{\nabla \eta_\e}^2 \to \SR^{N/2}$, $\ns{\eta_\e}^{2^*} \to \SR^{N/2}$ and
\be\lab{tem4-13-6}
	\f{\nm{\eta_\e}_{\al+\beta}^{\al+\beta}}{\nt{\eta_\e}^{\al+\beta-\p}}=\begin{cases}
	O(\e^{\f{6-\al-\beta}{4}})\quad&\text{if}~N=3,\\  O(|\log \e|^{-\f{\al+\beta-\p}{2}})\quad&\text{if}~N=4.
	\end{cases}
\ee
as $\e\to0$ from Lemma \ref{est} we get $e^{t_\e}\ge C\nt{\eta_\e}$. Substituting this into \eqref{tem4-13-8} we obtain
\be\lab{tem6-7-5}
	\begin{aligned}
	m_\nu(a,b)&\le\f{1}{N}\SR^{N/2}+O(\e^{N-2})+O(\e^{2\q})-
			C\e^{\beta\q}\f{\nm{\eta_\e}_{\al+\beta}^{\al+\beta}}{\nt{\eta_\e}^{\al+\beta-\p}}\\
		&\le\begin{cases} \f{1}{3}\SR^{3/2}+O(\e)+O(\e^{2\q}) -C\e^{\beta\q+\f{6-\al-\beta}{4}},\quad&\text{when}~N=3,\\
			\f{1}{4}\SR^2+O(\e^2) -C\e^{\beta}|\log \e|^{-\f{\al+\beta-\p}{2}},\quad&\text{when}~N=4.
		 \end{cases}
	\end{aligned}
\ee
This implies $m_\nu(a,b)\le\f{1}{N}\SR^{N/2}$ for all $\nu>0$ as claimed.\\

b) Using $\beta\q+\f{6-\al-\beta}{4}<\min\lbr{1,2\q}$ and $\beta<2$ it follows from \eqref{tem6-7-5} that $m_\nu(a,b)<\f{1}{N}\SR^{N/2}$ as required.\\

c) Arguing by by contradiction we assume $\nu_1=0$. Then there exists a sequence $\nu_n\to0^+$ such that
	$$0<m_{\nun}(a,b)<\f{1}{N}\SR^{N/2}\quad\text{for all $n>1$.}$$
Using Lemma \ref{PSseq} and Lemma \ref{property1}, we can apply Proposition \ref{PS2} with $c=m_\nun(a,b)$
to obtain a minimizer $(u_n,v_n)\in H_{rad}$ of $m_\nun(a,b)$ for any $n>1$.
It is easy to check that $(u_n,v_n)$ is also a positive solution of \eqref{mainequ}, that is:
\be\lab{tem6-7-6}
	\left\{ 	\begin{aligned}
	&-\dl u_n+\la_{1,n} u_n=|u_n|^{2^*-2}u_n+\al\nu_n |u_n|^{\al-2}|v_n|^\beta u_n,\quad \text{in }\RN,\\
	&-\dl v_n+\la_{2,n} v_n=|v_n|^{2^*-2}v_n+\beta\nu_n |u_n|^\al |v_n|^{\beta-2}v_n,\quad \text{in }\RN,\\
	\end{aligned}	\right.
\ee
with $\la_{1,n},\la_{2,n}>0$ and $\nt{u_n}=a$, $\nt{v_n}=b$.
Notice that
\be\lab{tem6-7-7}
	\begin{aligned}
		\f{1}{N}\SR^{N/2} &> m_\nun(a,b) = I_\nun(u_n,v_n)=I_\nun(u_n,v_n)-\f{1}{2}P_\nun(u_n,v_n)\\
			&=\f{1}{N}(\ns{u_n}^{2^*}+\ns{v_n}^{2^*})+\f{\p-2}{2}\nu_n\int_\RN|u_n|^\al|v_n|^\beta.
	\end{aligned}
\ee
Using the H\"older inequality we obtain 
\be
	\begin{aligned}
	\int_\RN |u_n|^\al|v_n|^\beta
		&\le \sbr{\int_\RN|u_n|^{2^*}}^{\f{\al}{2^*}} \sbr{\int_\RN|v_n|^{\f{ 2^*\beta}{2^*-\al }}}^{\f{2^*-\al}{2^*}}\\
		&\le \sbr{\int_\RN|u_n|^{2^*}}^{\f{\al}{2^*}} \sbr{\int_\RN|v_n|^2}^{\f{2^*-\al-\beta}{2^*-2}}
					\sbr{\int_\RN|v_n|^{2^*}}^{\f{1}{2^*} \f{2^*(\beta-2)+2\al}{2^*-2}}\\
		&\le \ns{u_n}^{\al} b^{\f{2^*-\al-\beta}{2^*-2}}
					\SR^{\f{N}{2}\f{1}{2^*} \f{2^*(\beta-2)+2\al}{2^*-2}}\\
		&=:Cb)\ns{u_n}^{\al}.
	\end{aligned}
\ee
Similarly we obtain $\int_\RN |u_n|^\al|v_n|^\beta \le Ca)\ns{v_n}^{\beta}$. Testing \eqref{tem6-7-6} with $(u_n,v_n)$ we deduce:
	$$\SR\ns{u_n}^2 \le \ns{u_n}^{2^*} + \nu_n\al Cb)\ns{u_n}^{\al}
	\quad\text{and}\quad
	\SR\ns{v_n}^2 \le \ns{v_n}^{2^*} + \nu_n\beta Ca)\ns{v_n}^{\beta}.$$
Since $\al,\beta\ge2$ a simple algebraic argument, see \cite[Lemma 3.3]{Peral etc=CVPDE=2009}, yields:
	$$\liminf_{n\to\infty}\ns{u_n}^{2^*} \ge \SR^{N/2} \quad\text{and}\quad \liminf_{n\to\infty}\ns{v_n}^{2^*} \ge \SR^{N/2}$$
contradicting \eqref{tem6-7-7} . Therefore $m_\nu(a,b)=\f{1}{N}\SR^{N/2}$ for $\nu>0$ small , i.e.\ $\nu_1>0$.\\

d) If $\nu_1>0$ then $m_\nu(a,b)=\f{1}{N}\SR^{N/2}$ for $0<\nu\le\nu_1$. If $m_\nu(a,b)$ is achieved by $(u,v)$ we can repeat the argument in b) and obtain a contradiction again.
\ep

\vskip0.12in
\bp[Proof of Theorem \ref{thm1} c) and d) ]
	The proof is just a combination of Proposition \ref{PS2}, Lemma \ref{PSseq}, Lemma \ref{property1} and Lemma \ref{property2}.
\ep

\vskip0.1in
\subsection{The case \texorpdfstring{$N=4$}{} and \texorpdfstring{$\al+\beta=3$}{}}

In this case we have $\ga_{\al+\beta}=2$. We define
\be\lab{nu0crit}
  \nu_0 := \f{N+2}{(\al+\beta)N\CR(N,\al,\beta) \sbr{a^2+b^2}^{2/N}}=\f{1}{2}\CR(4,\al,\beta)^{-1}(a^2+b^2)^{-1/2}	
\ee

\bl
	For $\nu>0$ we have 
\[
	\OR_\nu(a,b) := \lbr{(u,v)\in\TR(a,b): \nt{\nabla u}^2+\nt{\nabla v}^2>2\nu\int_\RN|u|^\al|v|^\beta} \neq \emptyset.
\] 
If $0<\nu<\nu_0$ then $\OR_\nu(a,b)=\TR(a,b)$.
\el

\bp
If $0<\nu<\nu_0$, then \eqref{VGNine} implies for $(u,v)\in\TR(a,b)$
\[
	2\nu\int_\RN|u|^\al|v|^\beta \le 2\nu\CR(4,\al,\beta)(a^2+b^2)^{1/2} (\nt{\nabla u}^2+\nt{\nabla v}^2)
		<\nt{\nabla u}^2+\nt{\nabla v}^2,
\]
hence $\TR(a,b) = \OR_\nu(a,b)$. If $\nu\ge\nu_0$ we claim that
\be\lab{tem6-11-1}
	\sup_{(u,v)\in\TR(a,b)}\f{\nt{\nabla u}^2+\nt{\nabla v}^2}{\int_\RN|u|^\al|v|^\beta}=+\iy,
\ee
from which $\OR_\nu(a,b)\neq\emptyset$ follows. Indeed, consider
	$$\phi_n(x):=\begin{cases} \sin n|x|\quad&\text{if}~|x|\le\pi\\
				0 \quad&\text{if}~|x|>\pi\end{cases}$$
and set $u_n:=\f{a}{\nt{\phi_n}}\phi_n$, $v_n:=\f{b}{\nt{\phi_n}}\phi_n$, so that $(u_n,v_n) \in \TR(a,b)$.
Since $\nt{\nabla \phi_n}^2= O(n^2)$, $\nt{\phi_n}^2=O(1)$ and $\np{\phi_n}^p=O(1)$, we have
	$$\f{\nt{\nabla u_n}^2+\nt{\nabla v_n}^2}{\int_\RN|u_n|^\al|v_n|^\beta}
			=\f{a^2+b^2}{a^\al b^\beta}\f{\nt{\phi_n}\nt{\nabla \phi_n}^2}{\np{\phi_n}^p}\to+\iy,\quad\text{as}~n\to\iy,$$
which implies \eqref{tem6-11-1}.
\ep

\bl\lab{structure3}
	For every $(u,v)\in \OR_\nu(a,b)$ the function $\Phi^\nu_{(u,v)}:\R\to\R$ has exactly one critical point $t_\nu(u,v)>0$.
	Moreover:
	\begin{itemize}
		\item[a)]	$\PR_\nu(a,b)=\PR_\nu^-(a,b)$ and $\PR_\nu(a,b)$ is a submanifold of $H$.
		\item[b)]	$t\s(u,v)\in\PR_\nu(a,b)$ if and only if $t=t_\nu(u,v)$.
		\item[c)] 	$\Phi^\nu_{(u,v)}(t)$ is strictly decreasing and concave on $(t_\nu(u,v),+\iy)$ and
						$$\Phi^\nu_{(u,v)}(t_\nu(u,v))=\max_{t\in\R}\Phi^\nu_{(u,v)}(t)>0.$$
		\item[d)] 	The map $(u,v)\mapsto t_\nu(u,v)$ is of class $C^1$.
	\end{itemize}
\el

\bl\lab{minimax2}
	$m_\nu(a,b)$ has a minimax  description:
		$$m_\nu(a,b)=\inf_{\OR_\nu(a,b)}\max_{t\in\R}I_\nu\big(t\s(u,v)\big).$$
\el

\bl\lab{PSseq2}
	If $0<\nu<\nu_0$ then $I_\nu|_{\TR(a,b)}$ has a Palais-Smale sequence $(u_n,v_n)\subset \TR(a,b)\cap H_{rad}$ at the
	level $m_\nu(a,b)$, with $P_\nu(u_n,v_n)\to0$ and $u_n^-,v_n^-\to0$ a.e. in $\RN$ as $n\to\iy$.
\el

	We omit the proofs of the lemmas \ref{structure3}--\ref{PSseq2} because they are similar to those of the lemmas \ref{structure2}--\ref{PSseq}.

\bl\lab{property3}
	Let $0<\nu<\nu_0$. For any fixed $a,b>0$, the following statements hold:
	\begin{itemize}
	\item[a)]	$m_\nu(a,b)\le m_\nu(a_1,b_1)$ for any $0< a_1\le a$ and $0<b_1\le b$;
	\item[b)]	$m_\nu(a,b)$ is nonincreasing with respect to $\nu\in(0,\nu_0)$;
	\end{itemize}
\el

\bp
a) Since $0<\nu<\nu_0$ implies $\OR_\nu(a,b)=\TR(a,b)$, the conclusion can be proved as in Lemma \ref{lem2-2}~b).\\

b) For any $0<\nu'\le\nu<\nu_0$, we see that
		$$m_\nu(a,b)=\inf_{\TR(a,b)}\max_{t\in\R}I_\nu\big(t\s(u,v)\big)
			\le \inf_{\TR(a,b)}\max_{t\in\R}I_{\nu'}\big(t\s(u,v)\big)=m_{\nu'}(a,b).$$
Therefore, $m_\nu(a,b)$ is nonincreasing with respect to $\nu\ge0$.
\ep

	We also have
\bl\lab{property4}
	Let $0<\nu<\nu_0$. For any fixed $a,b>0$,
		$$0<m_\nu(a,b)< \f{1}{4}\SR^{2}. $$
\el
\bp
For $(u,v)\in\PR_\nu(a,b)$ there holds:
	$$\begin{aligned}
	\nt{\nabla u}^2+\nt{\nabla v}^2&=\ns{u}^{2^*}+\ns{v}^{2^*}+2\nu \int_\RN|u|^\al|v|^\beta\\
		&\le \ns{u}^{2^*}+\ns{v}^{2^*}+\f{\nu}{\nu_0}(\nt{\nabla u}^2+\nt{\nabla v}^2)\\
		&\le C(\nt{\nabla u}^2+\nt{\nabla v}^2)^{\f{2^*}{2}}+\f{\nu}{\nu_0}(\nt{\nabla u}^2+\nt{\nabla v}^2) \\
	\end{aligned}$$
This implies $\inf_{\PR_\nu(a,b)}\sbr{\nt{\nabla u}^2+\nt{\nabla v}^2} >0$,
and therefore
	$$m_\nu(a,b)=\inf_{\PR_\nu(a,b)}\f{1}{N}(\ns{u_n}^{2^*}+\ns{v_n}^{2^*}) >0.$$

	Since $\al+\beta=3$, $\al>1$, $\beta>1$, we have $\al<2$, $\beta<2$. Recall the test functions $\eta_\e=\phi U_\e$ from Lemma~\ref{est} and set
$u_\e=\f{ a}{\nt{\eta_\e}}\eta_\e$ and $v_\e=\f{\e}{\nt{\eta_\e}}\eta_\e$.
Let $t_\e:=t_\nu(u_\e,v_\e)$ be given by Lemma \ref{structure3}. Then $t_\e\s(u_\e,v_\e)\in\PR_\nu(a,\e)$
and we obtain for $\e>0$ small:
\be \lab{tem6-12-2}
	\begin{aligned}
	m_\nu(a,b)&\le m_\nu(a,\e^\q)\le I_\nu(t_\e\s(u_\e,v_\e))\\
		&=\f{1}{2}e^{2t_\e}\sbr{\nt{\nabla u_\e}^2+\nt{\nabla v_\e}^2}
				-\f{1}{4}e^{4 t_\e}\sbr{\nm{u_\e}_4^{4}+\nm{v_\e}_4^{4}}
				-\nu e^{2t_\e}\int_\RN u_\e^\al v_\e^\beta\\
		&\le \max_{s>0}\sbr{ \f{a^2+\e^{2}}{2}s^2\f{\nt{\nabla\eta_\e}^2}{\nt{\eta_\e}^2}
			-\f{a^{4}+\e^{4}}{4}s^{4}\f{\nm{\eta_\e}_4^{4}}{\nt{\eta_\e}^{4}} }
			-C\e^{\beta} e^{2t_\e} \f{\nm{\eta_\e}_{3}^{3}}{\nt{\eta_\e}^{3}}\\
		&=\f{1}{N}\mbr{ \f{(a^2+\e^{2})\nt{\nabla\eta_\e}^2 }
			{\sbr{ a^{4}+\e^{4} }^{1/2}\nm{\eta_\e}_4^2 }}^{2}
			-C\e^{\beta} e^{2 t_\e} \f{\nm{\eta_\e}_{3}^{3}}{\nt{\eta_\e}^{3}}\\
		&=\f{1}{N}\SR^{N/2}+O(\e^2)
			-C\e^{\beta} e^{2t_\e} \f{\nm{\eta_\e}_3^3}{\nt{\eta_\e}^3}.
	\end{aligned}
\ee

	We claim that $\f{e^{t_\e}}{\nt{\eta_\e}}\ge C$ as $\e\to0$, for some constant $C>0$. By definition we have $P_\nu(t_\e\s(u_\e,v_\e))=0$, hence
\[
	\begin{aligned}
	e^{2t_\e}&=\f{\nt{\nabla u_\e}^2+\nt{\nabla v_\e}^2}{ \ns{u_\e}^{2^*}+\ns{v_\e}^{2^*}}
			-2\f{ \nu\int_\RN u_\e^\al v_\e^\beta }{\ns{u_\e}^{2^*}+\ns{v_\e}^{2^*}}\\
		&\ge \sbr{1-2\CR(N,\al,\beta)(a^2+b^2)^{1/2}\nu  }
			\f{\nt{\nabla u_\e}^2+\nt{\nabla v_\e}^2}{\ns{u_\e}^{2^*}+\ns{v_\e}^{2^*}}\\
		&\ge C \f{\nt{\nabla \eta_\e}^2\nt{\eta_\e}^{2^*-2}}{\ns{\eta_\e}^{2^*}},
	\end{aligned}
\]
where we used $0<\nu<\nu_0$. Now $\f{e^{t_\e}}{\nt{\eta_\e}}\ge C$ follows from Lemma \ref{est}.
Substituting $e^{t_\e}\ge C\nt{\eta_\e}$ into \eqref{tem6-12-2} we obtain
\be \lab{tem6-14-41}
	m_\nu(a,b) \le \f{1}{4}\SR^{2}+O(\e^2)-
			C\e^{\beta}\f{\nm{\eta_\e}_{3}^{3}}{\nt{\eta_\e}}
		\le\f{1}{4}\SR^2+O(\e^2) -C\e^{\beta}|\log \e|^{-1/2}.
\ee
From $\beta<2$ it follows that $m_\nu(a,b)<\f{1}{4}\SR^2$ for any $\e>0$ small, as claimed.
\ep

\br
	From \eqref{tem6-12-2} we see that $m_\nu(a,b)\le\f{1}{4}\SR^{2}f$ or any $\nu>0$.
\er

\bp[Proof of the Theorem \ref{thm1} b) and d)]
	The proof is just a combination of Proposition \ref{PS2}, Lemma \ref{PSseq2}, Lemma \ref{property3} and Lemma \ref{property4}.
\ep

\vskip0.3in
\section{Non-Existence in the repulsive case \texorpdfstring{$\nu<0$}{}}

We begin with a simple result about $\Phi_{(u,v)}^\nu:\R\to\R$.

\bl\lab{lem2-1}
	For $\nu\le0$ and $(u,v)\in \TR(a,b)$ there exists a unique $t_\nu(u,v) \in \R$ such that $t_\nu(u,v)\s(u,v) \in \PR_\nu(a,b)$;
	$t_\nu(u,v)$ is the unique critical point of $\Phi_{(u,v)}^\nu$ and a strict maximum point at positive level. Moreover:
	\begin{itemize}
		\item[a)]	$\PR_\nu(a,b)=\PR_\nu^-(a,b)$
		\item[b)]	$t_\nu(u,v)<0$ if and only if $P_\nu(u,v)<0$.
		\item[c)]	$\Phi_{(u,v)}^\nu$ is strictly decreasing and concave on $(t_\nu(u,v),+\iy)$.
		\item[d)]	The map $(u,v)\mapsto t_\nu(u,v)$ is of class $C^1$.
	\end{itemize}
\el
\bp
	Observe that
$$\sbr{\Phi_{(u,v)}^\nu}'(t)=e^{2^*t}\mbr{ \f{1}{e^{(2^*-2)t}}(\nt{\nabla u}^2+\nt{\nabla v}^2)
		+\f{|\nu|\p}{e^{(2^*-\p)t}} \int_\RN|u|^\al |v|^\beta-(\ns{u}^{2^*}+\ns{v}^{2^*}) }$$
is strictly decreasing because $\p<2^*$. It follows easily that $\sbr{\Phi_{(u,v)}^\nu}'$ has a unique zero $t_\nu(u,v)$, hence $t_\nu(u,v)\s (u,v)\in\PR_\nu(a,b)$. Moreover, $\sbr{\Phi_{(u,v)}^\nu}'(t)>0$ for $t < t_\nu(u,v)$, and that $\sbr{\Phi_{(u,v)}^\nu}'(t)<0$ for $t > t_\nu(u,v)$, so
$t_\nu(u,v)$ is a strict maximum point of $\Phi_{(u,v)}^\nu$ on positive level. The conclusions $a)$-$d)$ easily follow from the properties of $\Phi^\nu_{(u,v)}$.
\ep

Setting
	$$\BR_\nu(a,b):=\lbr{ (u,v)\in \TR(a,b):P_\nu(u,v)\le0 },$$
	$$\UR_\nu(a,b):=\lbr{ (u,v)\in \TR(a,b):P_\nu(u,v)<0 },$$
and
	$$J_\nu(u,v) := I_\nu(u,v)-\f{1}{2^*}P_\nu(u,v) = \f{1}{N}(\nm{\nabla u}^2+\nm{\nabla v}^2)-\f{2^*-\p}{2^*}\nu\int_\RN|u|^\al|v|^\beta.$$
Then we have the following lemma.
\bl\lab{lem2-2}
	Let $\nu\le0$ and $a,b>0$. Then the following holds.
	\begin{itemize}
	\item[a)]	$m_\nu(a,b)=\inf_{\BR_\nu(a,b)} J_\nu(u,v)=\inf_{\UR_\nu(a,b)} J_\nu(u,v)$.
	\item[b)]	If $0<a_1\le a$, $0<b_1\le b$ then $m_\nu(a,b)\le m_\nu(a_1,b_1)$.
	\item[c)]	For each $\nu\le0$ there holds $m_\nu(a,b)= \f{1}{N}\SR^{\f{N}{2}}$ with $\SR$ being defined in \eqref{Sobolev1}.
	\end{itemize}
\el

\bp
a) Observe that $\UR_\nu(a,b)$ is dense in $\BR_\nu(a,b)$ by Lemma~\ref{lem2-1}, hence 
\[
  \inf_{\UR_\nu(a,b)}J_\nu = \inf_{\BR_\nu(a,b)}J_\nu \le \inf_{\PR_\nu(a,b)} J_\nu = \inf_{\PR_\nu(a,b)}I_\nu = m_\nu(a,b)
\]
where we also used that $\PR_\nu(a,b) \subset \BR_\nu(a,b)$ and $P_\nu=0$ on $\PR_\nu(a,b)$. Moreover, for $(u,v)\in\UR_\nu(a,b)$ we have $t_\nu(u,v)<0$ according to Lemma \ref{lem2-1}. This implies for any $(u,v)\in\UR_\nu(a,b)$
\begin{align*} 
	J_\nu(u,v)& > \f{1}{N}e^{2t_\nu(u,v)}(\nt{\nabla u}^2+\nt{\nabla v}^2)
			       +\f{2^*-\p}{2^*}|\nu|e^{\p t_\nu(u,v)}\int_\RN|u|^\al|v|^\beta\\
		&= J_\nu(t_\nu(u,v)\s (u,v)) \ge \inf_{\PR_\nu(a,b)}J_\nu
\end{align*}
hence $\displaystyle \inf_{\UR_\nu(a,b)}J_\nu \ge \inf_{\PR_\nu(a,b)}J_\nu$.\\

b) It is sufficient to prove that for $0<a_1\le a$, $0<b_1\le b$ and arbitrary $\e>0$, one has
\be\lab{tem2-4}
	m_\nu(a,b)\le m_\nu(a_1,b_1)+\e.
\ee
By the definition of $m_\nu(a_1,b_1)$ there exists $(u,v)\in\PR_\nu(a_1,b_1)$ such that
\be\lab{tem2-1}
	I_\nu(u,v)\le m_\nu(a_1,b_1)+\f{\e}{2}.
\ee
Let $\vp\in\CR_0^\iy(\RN)$ be radial and such that
	$$0\le\phi\le1\quad\text{and}\quad \phi(x)=\begin{cases}1\quad\text{if}~|x|\le1,\\0\quad\text{if}~|x|\ge2. \end{cases}$$
For $\delta>0$ we define $u_\delta(x):=u(x)\vp(\delta x)$ and $v_\delta(x):=v(x)\vp(\delta x)$ so that $(u_\delta,v_\delta) \to (u,v)$ 
in $H$ as $\delta\to0^+$. Then Lemma \ref{lem2-1}~d) implies $t_\nu(u_\delta,v_\delta)\s(u_\delta,v_\delta)\to t_\nu(u,v)\s(u,v)$  in $H$ as $\delta\to0^+$.
As a consequence, we can fix $\delta>0$ such that
\be\lab{tem2-2}
	I_\nu(t_\nu(u_\delta,v_\delta)\s(u_\delta,v_\delta))\le I_\nu(u,v)+\f{\e}{4}.
\ee
Next we choose $\psi\in\CR_0^\iy(\RN)$ with $supp(\psi)\subset \R^N\setminus B(0,\f{4}{\delta})$ and set
	$$\psi_a=\f{\sqrt{a^2-\nt{u_\delta}^2}}{\nt{\psi}}\psi\quad\text{and}
		\quad\psi_b=\f{\sqrt{b^2-\nt{v_\delta}^2}}{\nt{\psi}}\psi.$$
Then
	$$\sbr{supp(u_\delta)\cup supp(u_\delta)}\cap\sbr{supp(\tau\s\psi_a)\cup supp(\tau\s\psi_b)}=\emptyset \quad\text{for $\tau\le0$,}$$
and thus
	$$(\tilde u_\tau,\tilde v_\tau):=(u_\delta+\tau\s\psi_a,v_\delta+\tau\s\psi_b)\in \TR(a,b).$$
Setting $t_\tau:=t_\nu(\tilde u_\tau,\tilde v_\tau)$ we have $P_\nu(t_\tau\s(\tilde u_\tau,\tilde v_\tau))=0$, i.e.
	$$\f{1}{e^{(2^*-2)t_\tau}}\Big(\nt{\nabla \tilde u_\tau}^2+\nt{\nabla \tilde v_\tau}^2\Big)
		+\f{|\nu|\p}{e^{(2^*-\p)t_\tau}}\int_\RN|\tilde u_\tau|^\al|\tilde v_\tau|^\beta
		=\ns{\tilde u_\tau}^{2^*}+\ns{\tilde v_\tau}^{2^*}.$$
It follows that $\limsup_{\tau\to-\iy} t_\tau<+\iy$, because $(\tilde u_\tau,\tilde v_\tau)\to(u_\delta,v_\delta)\neq(0,0)$ as
$\tau\to-\iy$. Consequently $t_\tau+\tau\to-\iy$ as $\tau\to-\iy$ and we have for $\tau<-1$ small:
\be\lab{tem2-3}
	I_\nu\big((t_\tau+\tau)\s(\psi_a,\psi_b)\big)<\f{\e}{4}.
\ee
Now using \eqref{tem2-1}, \eqref{tem2-2} and \eqref{tem2-3}, we obtain
\begin{align*}
	m_\nu(a,b)&\le I_\nu\big(t_\tau\s(\tilde u_\tau,\tilde v_\tau)\big)
		 =I_\nu\big(t_\tau\s(u_\delta,v_\delta)\big)+I_\nu\big((t_\tau+\tau)\s(\psi_a,\psi_b)\big)\\
		&\le I_\nu\big(t_\nu(u_\delta,v_\delta)\s(u_\delta,v_\delta)\big)+\f{\e}{4}
		 \le I_\nu(u,v)+\f{\e}{4}\le m_\nu(a_1,b_1)+\e,
\end{align*}
that is \eqref{tem2-4}.\\

c) We prove the conclusion in two steps.
\begin{itemize}[fullwidth,itemindent=1em]
\vskip 0.05in
\item[\bf{Step 1.)}] We prove that $m_\nu(a,b)=m_0(a,b)$ for $\nu<0$.\\
	Since $P_0(u,v)\le P_\nu(u,v)$ and $J_0(u,v)\le J_\nu(u,v)$  we have $\BR_\nu(a,b)\subset\BR_0(a,b)$ and
	$$m_0(a,b)=\inf_{\BR_0(a,b)}J_0(u,v)\le \inf_{\BR_\nu(a,b)}J_\nu(u,v)=m_\nu(a,b).$$
On the other hand, for $(u,v)\in\UR_0(a,b)$ and $t>0$ we set $\big(u_t(x),v_t(x)\big)=\big(t^{\f{N-2}{2}}u(tx),t^{\f{N-2}{2}}v(tx)\big)$
and obtain for $t>1$ karge:
\begin{align*}
	P_\nu(u_t,v_t)&=\nt{\nabla u}^2+\nt{\nabla v}^2-\ns{u}^{2^*}-\ns{v}^{2^*} + \f{\p|\nu|}{t^{\al+\beta-\p}}\int_\RN|u|^\al|v|^\beta\\
			&=P_0(u,v)+\f{\p|\nu|}{t^{\al+\beta-\p}}\int_\RN|u|^\al|v|^\beta<0.
\end{align*}
Noting that $\nt{u_t}=t^{-1}a$, $\nt{v_t}=t^{-1}b$, we get for $t>1$:
\begin{align*}
	m_\nu(a,b)&\le m_\nu(t^{-1}a,t^{-1}b)\le J_\nu(u_t,v_t)
		=\f{1}{N}(\nt{\nabla u_t}^2+\nt{v_t}^2) + \f{2^*-\p}{2^*}|\nu|\int_\RN|u_t|^\al|v_t|^\beta\\
		&=J_0(u,v)+\f{2^*-\p}{2^*}\f{|\nu|}{t^{\al+\beta-\p}}\int_\RN|u|^\al|v|^\beta
		\to J_0(u,v),\quad \text{as}~t\to\iy.
\end{align*}This implies $m_\nu(a,b)\le m_0(a,b)$ as claimed.

\vskip 0.05in
\item[\bf{Step 2.)}] We prove that $m_0(a,b)=\f{1}{N}\SR^{\f{N}{2}}$. \\
For $(u,v)\in \PR_0(a,b)$, we have
	$$ \nt{\nabla u}^2+\nt{\nabla v}^2
	    = \ns{u}^{2^*}+\ns{v}^{2^*}
	    \le \SR^{-\f{2^*}{2}}\big(\nt{\nabla u}^{2^*} +\nt{\nabla v}^{2^*}\big)
	   < \SR^{-\f{2^*}{2}}\big(\nt{\nabla u}^{2} + \nt{\nabla v}^{2}\big)^{\f{2^*}2}, $$
and therefore $m_0(a,b)\ge \f{1}{N}\SR^{\f{N}{2}}$.
On the other hand, according to \cite[Proposition 2.2]{Soave=JFA=2020}, given $\e>0$ there exists $u\in H^1(\RN)$ with $\nt{u}=a$ such that 
\be\lab{tem2-6}
	\max_{t\in\R}E(t\s u)\le \f{1}{N}\SR^{\f{N}{2}}+\e,
\ee
where $E(u)$ is defined in \eqref{defofE}. Now choose any $v\in H^1(\RN)$ with $\nt{v}=b$, define $t(s):=t_0(u,s\s v)$, and observe that $\limsup_{\tau\to-\iy} t(s)<+\iy$ as in b). We obtain:
\begin{align*}
	m_0(a,b)&\le I_0\big(t(s)\s(u,s\s v))
		       = E\big(t(s)\s u\big) + E\big((t(s)+s)\s v\big)\\
		&\le \f{1}{N}\SR^{\f{N}{2}}+\e + E\big((t(s)+s)\s v\big) \to \f{1}{N}\SR^{\f{N}{2}}+\e \quad\text{as $s\to-\iy$.}
\end{align*}
Since this holds for any $\e>0$ we deduce $m_0(a,b)\le\f{1}{N}\SR^{\f{N}{2}}$.
\end{itemize}
\ep

\bp[Proof of Theorem \ref{thm2}.]
For $(u,v)\in\PR_\nu(a,b)$ we have for $\nu\le0$:
\[
	\begin{aligned}
	 \nt{\nabla u}^2+\nt{\nabla v}^2&=\ns{u}^{2^*}+\ns{v}^{2^*}+\p \nu\int_\RN|u|^\al|v|^\beta \le \ns{u}^{2^*}+\ns{v}^{2^*}\\
	 	&\le \SR^{-\f{2^*}{2}}\big(\nt{\nabla u}^{2^*} + \nt{\nabla v}^{2^*}\big)	
	 	  < \SR^{-\f{2^*}{2}}\big(\nt{\nabla u}^{2} + \nt{\nabla v}^{2}\big)^{\f{2^*}2}
	\end{aligned}
\]
so that $\nt{\nabla u}^2+\nt{\nabla v}^2 > \f{1}{N}\SR^{\f{N}{2}}$. This implies for $\nu \le 0$:
\begin{align*}
	 I_\nu(u,v) &= I_\nu(u,v)-\f1{2^*}\PR_\nu(u,v)
	                  = \f1N\big(\nt{\nabla u}^{2} + \nt{\nabla v}\big) - \f{2^*-\p}{2}\nu\int_{\RN}|u|^\al|v|^\beta\\
	                &\ge \f1N\big(\nt{\nabla u}^{2} + \nt{\nabla v}^2\big) 
	                  > \f{1}{N}\SR^{\f{N}{2}} = m_\nu(a,b)
\end{align*}
\ep

\vskip0.3in
\section{The limit system}

Inspired by \cite{Szulkin-Weth=JFA=2009} we consider for $(u,v)\in\TR(a,b)$ the equation
\[
  \f{d}{dt} K\big(t\s(u,v)\big) = \f{d}{dt} \left(\f{e^{2t}}{2}\big(\nt{\nabla u}^2+\nt{\nabla v}^2\big) - e^{\p t}\int_{\RN}|u|^\al|v|^\beta\right) = 0
\]
which has the solution $t(u,v)$ defined by
\[
  e^{t(u,v)} = \sbr{\f{\nt{\nabla u}^2+\nt{\nabla v}^2}{\p\nm{|u|^\al|v|^\beta}_1}}^{\f{1}{2-\p}}.
\]
Then $t(u,v)\s(u,v)\in\LR(a,b)$ by definition \eqref{limitpho} and
\be\lab{9290}
	\tilde K(u,v) := K\big(t(u,v)\s(u,v)\big) = \f{\p-2}{2\p}
			\mbr{ \f{\sbr{\nt{\nabla u}^2+\nt{\nabla v}^2}^\f{\p}{2}}{\p\nm{|u|^\al|v|^\beta}_1 }}^{\f{2}{\p-2}}.
\ee
Moreover $l(a,b)=\inf_{\LR(a,b)}K$ satisfies the following.

\bl\lab{limit1} $l(a,b) = \inf_{\TR(a,b)}\tilde K$ for all $a,b>0$.
\el

\bp
	For $(u,v)\in\TR(a,b)$ we have $\tilde K(u,v) \ge \inf_{\LR(a,b)} K = l(a,b)$. On the other hand
	$$\inf_{\TR(a,b)}\tilde K(u,v)\le \inf_{\TR(a,b)}K(u,v)\le l(a,b),$$
where in the last inequality we use $\LR(a,b)\subset\TR(a,b)$.
\ep

\bl\lab{limit2}
	\begin{itemize}
	\item[a)]	If $\al+\beta<2+\f4N$ then $l(a,b)<0$.
	\item[b)]	If $\al+\beta>2+\f4N$ then $l(a,b)>0$.
	\item[c)]	$l(a,b)\le l(a_1,b_1)$ for $0< a_1\le a$ and $0<b_1\le b$.
	\end{itemize}
\el

\bp
a) If $\al+\beta<2+\f4N$ then $\p<2$, hence $\tilde K(u,v)<0$ for any $(u,v)\in \TR(a,b)$.\\

b) Using \eqref{VGNine}, we have for $(u,v)\in\LR(a,b)$:
	$$\nt{\nabla u}^2+\nt{\nabla v}^2=\int_\RN|u|^\al|v|^\beta\le C\Big(\nt{\nabla u}^2+\nt{\nabla v}^2\Big)^{\p/2},$$
which implies $\inf_{\LR(a,b)}\nt{\nabla u}^2+\nt{\nabla v}^2>0$, and therefore
	$$l(a,b)=\inf_{\LR(a,b)}\sbr{\f{1}{2}-\f{1}{\p}}(\nt{\nabla u}^2+\nt{\nabla v}^2)>0.$$

c) If $\al+\beta<2+\f4N$ we obtain as in \cite[Lemma 3.1]{Gou-Jeanjean=2016} for $0 < a_1< a$ and $0<b_1 < b$:
	$$l(a,b)\le l(a_1,b_1)+l\Big(\sqrt{a^2-a_1^2}, \sqrt{b^2-b_1^2}\Big) < l(a_1,b_1).$$
	because $l\Big(\sqrt{a^2-a_1^2}, \sqrt{b^2-b_1^2}\Big)<0$ by a).
	In the case $\al+\beta>2+\f4N$ we can argue as in the proof of Lemma \ref{lem2-2} (2).
\ep

\bp[Proof of Theorem \ref{thmlim}]
If $\al+\beta<2+\f4N$ we obtain a Palais-Smale sequence by starting with a minimizing sequence $K(u_n,v_n)\to l(a,b)$ with $(u_n,v_n)\in\TR(a,b)$. Using symmetric decreasing rearrangement we may assume that $u_n$ and $v_n$ are non-negative and radially symmetric. By Ekeland's variational principle we may also assume that $(u_n,v_n)$ is a Palais-Smale sequence. Then we have
	$$K'(u_n,v_n)+\la_{1,n}u_n+\la_{2,n}v_n\to0$$
with
	$$\la_{1,n}=-\f{1}{a^2}K'(u_n,v_n)[(u_n,0)]\quad\text{and}\quad
		\la_{2,n}=-\f{1}{b^2}K'(u_n,v_n)[(0,v_n)]. $$
Using \eqref{VGNine}, we see that $K$ is coercive on $\TR(a,b)$, hence $(u_n,v_n)$ is bounded in $H$ and $\la_{1,n},\la_{2,n}$ are bounded in $\R$.
Then there exist $(u,v)\in H_{rad}$ and $\la_1,\la_2\in\R$ such that up to a subsequence:
	$$(u_n,v_n)\rh(u,v)\quad\text{in} ~H_{rad},$$
	$$(u_n,v_n)\ra(u,v)\quad\text{in} ~L^q(\RN)\times L^q(\RN),~\text{for}~2<q<2^*,$$
	$$(u_n,v_n)\ra(u,v)\quad\text{a.e. in}~\RN,$$
	$$(\la_{1,n},\la_{2,n})\to(\la_1,\la_2)\quad\text{in} ~\R^2.$$
We claim that $u\neq0$ and $v\neq0$. Indeed, if $u=0$ or $v=0$, then
	$$l(a,b)=\lim_{n\to\iy}\f{1}{2}(\nt{\nabla u_n}^2+\nt{\nabla v_n}^2)-\nm{|u_n|^\al|v_n|^\beta}_1\ge0,$$
a contradiction with $l(a,b)<0$. So we have $0<\nt{u}\le a$ and $0<\nt{v}\le b$.
Moreover, $(u,v)$ is a solution of
\be
	\begin{cases}
	&-\Delta u+\la_1u=\al|u|^{\al-2}|v|^\beta u,\quad\text{in}~\RN,\\
	&-\Delta v+\la_2v=\beta|u|^\al|v|^{\beta-2} v,\quad\text{in}~\RN,\\
	&u\ge0, \quad v\ge0.
	\end{cases}
\ee
We claim that $\la_1,\la_2>0$. Indeed, if $\la_1\le0$ then
	$$-\Delta u=|\la_1|u+\al|u|^{\al-2}|v|^\beta u\ge 0,\quad\text{in}~\RN,$$
and hence $u=0$ by \cite[Lemma A.2]{Ikoma}, which is impossible. Analogously we can prove $\la_2>0$.
Note that $\nabla K(u,v)+(\la_1 u,\la_2 v)=0$, so
	$$\begin{aligned}
		&\nt{\nabla (u_n-u)}^2+\la_1\nt{u_n-u}^2\\
		&\hspace{1cm}= \mbr{K'(u_n,v_n)-K'(u,v)}[(u_n-v,0)]+\la_{1,n}\int_\RN u_n(u_n-u)-\la_1 \int_\RN u(u_n-u)+o_n(1)\\
		&\hspace{1cm}= o_n(1),	\end{aligned}$$
which implies $u_n\to u$ in $H^1(\RN)$, using $\la_1>0$. Analogously we can prove $v_n\to v$ in $H^1(\RN)$. 
Then $(u,v)\in\LR(a,b)$, and $K(u,v)=l(a,b)$. Clearly, $(u,v)$ is a solution of \eqref{limitequ} and therefore a normalized ground state of \eqref{limitequ}.
Moreover, $u,v$ are positive by the maximum principle, and they are radially symmetric.

If $\al+\beta>2+\f4N$ one can construct a Palais-Smale sequence $(u_n,v_n)\in\LR(a,b)$ for $K|_{\TR(a,b)}$ at the level $l(a,b)>0$ with $L(u_n,v_n)\to0$
and $u_n^-,v_n^-\to0$ a.e. in $\RN$ as $n\to\iy$. Then one can argue as in the case $\al+\beta<2+\f4N$.
\ep

\bp[Proof of Proposition \ref{prop2}]
	Let $(u,v)\in\LR(a,b)$ be such that $K(u,v)=l(a,b)$. Then
		$$l(a,b)=\f{\p-2}{2\p}(\nt{\nabla u}^2+\nt{\nabla v}^2)=\f{\p-2}{2}\nm{|u|^\al|v|^\beta}_1,$$
and
	$$\tilde K(u,v)=\f{\p-2}{2\p}(\nt{\nabla u}^2+\nt{\nabla v}^2)=l(a,b).$$
Due to definition \eqref{VGNconstant}, we also have
\be\lab{limit-3}
	\begin{aligned}
	\CR(N,\al,\beta)^{-1}&\le Q(u,v)=\p(a^2+b^2)^{(\al+\beta-\p)/2}(\nt{\nabla u}^2+\nt{\nabla v}^2)^{\f{\p-2}{2}}\\
		&=\p(a^2+b^2)^{(\al+\beta-\p)/2}\sbr{\f{2\p}{\p-2}l(a,b)}^{\f{\p-2}{2}}.
	\end{aligned}
\ee
On the other hand, since $\f{a^2}{b^2}=\f{\al}{\beta}$, there exists $\sigma_1>0$ such that
$(Z_1^{\sigma_1,1},Z_2^{\sigma_1,1})\in\TR(a,b)$, where $(Z_1,Z_2)$ is defined in Proposition \ref{VGN}.
It follows that
\be\lab{limit-4}
	l(a,b)\le \tilde K(Z_1^{\sigma_1,1},Z_2^{\sigma_1,1})
		=\f{\p-2}{2\p}\p^{\f{2}{2-\p}}(a^2+b^2)^{\f{(\al+\beta-\p)}{2-\p}}\CR(N,\al,\beta)^{\f{2}{2-\p}}.
\ee
Combining \eqref{limit-3} and \eqref{limit-4} we obtain
\be\lab{limit-5}
	l(a,b)=\f{\p-2}{2\p}\p^{\f{2}{2-\p}}(a^2+b^2)^{\f{(\al+\beta-\p)}{2-\p}}\CR(N,\al,\beta)^{\f{2}{2-\p}}.
\ee
Now a normalized ground state $(W_1,W_2)\in\GR(a,b)$ of \eqref{limitequ} must satisfy $(W_1,W_2)\in\LR(a,b)$.
Since $l(a,b)$ is achieved by a normalized ground state of \eqref{limitequ}, we have $K(W_1,W_2)=l(a,b)$.
Then from \eqref{limit-4} and \eqref{limit-5}, we conclude that $Q(W_1,W_2)=\CR(N,\al,\beta)^{-1}$.
According to Proposition \ref{VGN} we have
	$$(W_1,W_2)=(Z_1^{\sigma,\mu}, Z_2^{\sigma,\mu})\quad\text{for some}~\sigma>0, \mu>0.$$
Substituting $(u_0,v_0)$ into the equation \eqref{limitequ} we get
\[
	\begin{cases}
	\sigma^2\mu^{-N}\al^{-\f{2-\beta}{\al+\beta-2}}\beta^{-\f{\beta}{\al+\beta-2}}\nt{Z}^2=a^2\\
	\sigma^2\mu^{-N}\al^{-\f{\al}{\al+\beta-2}}\beta^{-\f{2-\al}{\al+\beta-2}}\nt{Z}^2=b^2,\\
	\sigma^{2-\al-\beta}\mu^2=\f{1}{\al+\beta}.
	\end{cases}
\]
Since $\f{a^2}{b^2}=\f{\al}{\beta}$ this has the unique solution
\be\lab{mu}
	\begin{cases}
		\sigma&=\sbr{\f{a}{\nt{Z}}}^{\f{2}{2-\p}}\mbr{ \sbr{\f{\al}{\beta}}^{\beta/2}+\sbr{\f{\beta}{\al}}^{1-\beta/2} }
					^{\f{1}{\p-2}\f{2}{\al+\beta-2}}(\al+\beta)^{\f{1}{\al+\beta-2}},\\
		\mu&=\sbr{\f{a}{\nt{Z}}}^{\f{\al+\beta-2}{2-\p}}\mbr{ \sbr{\f{\al}{\beta}}^{\beta/2}+\sbr{\f{\beta}{\al}}^{1-\beta/2} }
					^{\f{1}{\p-2}}.	
	\end{cases}
\ee
It follows that $W_i(x)=\sigma_i Z(\mu x)$ for $i=1,2$ with
\be\lab{sigma}
	\sigma_1=\sigma \al^{-\f{2-\beta}{2(\al+\beta-2)}} \beta^{-\f{\beta}{2(\al+\beta-2)}}\quad\text{and}\quad
	\sigma_2=\sigma \al^{-\f{\al}{2-(\al+\beta-2)}} \beta^{-\f{2-\beta}{2(\al+\beta-2)}}.
\ee
\ep

\vskip0.3in
\section{The asymptotic behavior}

	In the setting of Theorem \ref{thm1} we know that $m_\nu(a,b)$ is achieved by some $(u_\nu,v_\nu)$ for $\nu>0$ small or $\nu>0$ large,
which is positive, radially symmetric and radially decreasing. Our goal in this section is to study the asymptotic behavior of $(u_\nu,v_\nu)$ as $\nu\to0^+$ or $\nu\to+\iy$, and to prov Theorems~\ref{thm3} and \ref{thm4}.
 We recall that $(u_\nu,v_\nu)$ is also a normalized
ground state of \eqref{mainequ} with $\la_{1,\nu},\la_{2,\nu}>0$.
For simplicity, we write $A_\nu \sim B_\nu$, $A_\nu \lesssim B_\nu$ and $A_\nu \gtrsim B_\nu$ if there exist $C_1,C_2>0$ such that, respectively, $C_1 B_\nu \le A_\nu \le C_2 B_\nu$, $A_\nu \le C_1 B_\nu$ and $A_\nu \ge C_1B_\nu$ as $\nu \to 0+$ (or $\nu \to \iy$).

\bl\lab{tem1}
	If $\al+\beta<2+\f{4}{N}$ then
		$$\la_{1,\nu}+\la_{2,\nu}\sim \nt{\nabla u_\nu}^2+\nt{\nabla v_\nu}^2\sim \nu^{\f{2}{2-\p}},$$
	as $\nu\to0^+$.
\el
\bp
	Since $(u_\nu,v_\nu)\in\PR_\nu^+(a,b)$, we have
\be\lab{tem6-14-1}
	\nt{\nabla u_\nu}^2+\nt{\nabla v_\nu}^2 = \ns{u_\nu}^{2^*}+\ns{v_\nu}^{2^*}+\p\nu\nm{u_\nu^\al v_\nu^\beta}_1,
\ee
and
\be\lab{tem6-14-4}
	2\big(\nt{\nabla u_\nu}^2+\nt{\nabla v_\nu}^2\big) > 2^*(\ns{u_\nu}^{2^*}+\ns{v_\nu}^{2^*})
				+\p^2\nu\nm{u_\nu^\al v_\nu^\beta}_1.
\ee
It follows from \eqref{VGNine} that
	$$\nt{\nabla u_\nu}^2+\nt{\nabla v_\nu}^2\lesssim \nu\nm{u_\nu^\al v_\nu^\beta}_1
			\lesssim \nu\big(\nt{\nabla u_\nu}^2+\nt{\nabla v_\nu}^2\big)^{\f{\p}{2}},$$
which implies
\be\lab{tem6-14-2}
	\nt{\nabla u_\nu}^2+\nt{\nabla v_\nu}^2\lesssim \nu^{\f{2}{2-\p}}
\ee
and, by \eqref{tem6-14-1} 
\be\lab{tem6-14-3}
	\nu\nm{|u_\nu|^\al |v_\nu|^\beta}_1\lesssim \nu^{\f{2}{2-\p}}.
\ee
For $(u,v)\in\TR(a,b)$ and $\nu>0$ small there exists by Lemma \ref{structure1} a unique $t_\nu:=t_\nu(u,v)$ such that
 $t_\nu\s(u,v)\in\PR_\nu^+(a,b)$, i.e. 
	$$e^{2t_\nu}(\nt{\nabla u}^2+\nt{\nabla v}^2)=e^{2^*t_\nu}(\ns{u}^{2^*}+\ns{v}^{2^*})
			+e^{\p t_\nu}\p\nu\nm{|u|^\al |v|^\beta}_1,$$
and
	$$2e^{2t_\nu}(\nt{\nabla u}^2+\nt{\nabla v}^2)> 2^*e^{2^*t_\nu}(\ns{u}^{2^*}+\ns{v}^{2^*})
				+e^{\p t_\nu}\p^2\nu\nm{|u|^\al |v|^\beta}_1.$$
Therefore
	$$(2^*-2)e^{2t_\nu}(\nt{\nabla u}^2+\nt{\nabla v}^2)< \nu(2^*-\p)\p e^{\p t_\nu}\nm{|u|^\al |v|^\beta}_1,$$
which implies $e^{t_\nu}\to0$ as $\nu\to0$ because $\p<2$. It follows that
	$$e^{2t_\nu}\sim \nu e^{\p t_\nu},\quad\text{as}~\nu\to0,$$
which implies $e^{t_\nu}\sim \nu^{\f{1}{2-\p}}$ as $\nu\to0$, hence
	$$\begin{aligned}
		I_\nu\big(t_\nu\s(u,v)\big)&=\f{\p-2}{2\p}e^{2t_\nu}\big(\nt{\nabla u}^2+\nt{\nabla v}^2\big)
				+\f{2^*-\p}{2^*\p}e^{2^*t_\nu}\big(\ns{u}^{2^*}+\ns{v}^{2^*}\big)\\
			&\sim -\nu^{\f{2}{2-\p}}.
	\end{aligned}$$
Using $I_\nu(t_\nu\s(u,v))\ge m_\nu(a,b)$ and $m_\nu(a,b)\gtrsim -\nu\nm{u_\nu^\al v_\nu^\beta}_1$ we next deduce
	$$\nu\nm{u_\nu^\al v_\nu^\beta}_1\gtrsim \nu^{\f{2}{2-\p}},$$
which together with \eqref{tem6-14-3}, implies
	$$\nu\nm{u_\nu^\al v_\nu^\beta}_1\sim \nu^{\f{2}{2-\p}}.$$
The Pohozaev identity yields
	$$\la_{1,\nu}+\la_{2,\nu}\sim \nu\nm{u_\nu^\al v_\nu^\beta}_1\sim \nu^{\f{2}{2-\p}},$$
and by \eqref{tem6-14-4} and \eqref{tem6-14-2} we get
	$$\nt{\nabla u_\nu}^2+\nt{\nabla v_\nu}^2\sim \nu\nm{u_\nu^\al v_\nu^\beta}_1\sim \nu^{\f{2}{2-\p}}$$
which completes the proof.
\ep

	By ??? we see that for any $(u,v)\in\TR(a,b)$,
\be\lab{est1}
	\nm{|u|^\al |v|^\beta}_1\le D_0 (\nt{\nabla u}^2+\nt{\nabla v}^2)^{\f{\p}{2}}
\ee
with
	$$D_0:=\p^{-1}\sbr{\f{\p-2}{2\p l(a,b)}}^{\f{\p-2}{2}}>0.$$

\bl\lab{tem4}
	If $\al+\beta<2+\f{4}{N}$ then
		$$dist_H\sbr{(-t_\nu)\s(u_\nu,v_\nu),\GR(a,b)}\to0\quad\text{as}~\nu\to0^+,$$
	where $t_\nu:=t_\nu(u_0,v_0)$ is from Lemma \ref{structure1}. Moreover we have $e^{t_\nu}\sim \nu^{\f{1}{2-\p}}$ as $\nu\to0^+$.
\el

\bp
	Let $(u_0,v_0)\in\GR(a,b)$. Then $\tilde K(u_0,v_0)=l(a,b)$ which implies
\be\lab{tem6-14-13}
	\nm{u_0^\al v_0^\beta}_1=D_0 (\nt{\nabla u_0}^2+\nt{\nabla v_0}^2)^{\f{\p}{2}}.
\ee
Since $t_\nu \in \PR_\nu^+(a,b)$ we have
\be\lab{tem6-14-23}
	e^{2t_\nu}\big(\nt{\nabla u_0}^2+\nt{\nabla v_0}^2\big)=e^{2^*t_\nu}\big(\ns{u_0}^{2^*}+\ns{v_0}^{2^*}\big)
			+e^{\p s_\nu}\p\nu\nm{u_0^\al v_0^\beta}_1,
\ee
and
\be\lab{tem6-14-24}
	2e^{2s_\nu}(\nt{\nabla u_0}^2+\nt{\nabla v_0}^2)> 2^*e^{2^*s_\nu}(\ns{u_0}^{2^*}+\ns{v_0}^{2^*})
				+e^{\p s_\nu}\p^2\nu\nm{u_0^\al v_0^\beta}_1.
\ee
Using \eqref{est1} we have
\be\lab{tem6-14-16}
	\begin{aligned}
		&(2^*-2)(\nt{\nabla (t_\nu\s u_0)}^2+\nt{\nabla (t_\nu\s v_0)}^2)
		  < (2^*-\p)\p\nu\nm{|t_\nu\s u_0|^\al |t_\nu\s v_0|^\beta}_1\\
		&\hspace{2cm}\le (2^*-\p)\p\nu D_0(\nt{\nabla (t_\nu\s u_0)}^2+\nt{\nabla (t_\nu\s v_0)}^2)^{\f{\p}{2}},
	\end{aligned}
\ee
which yields
\be\lab{tem6-14-11}
	\nt{\nabla (t_\nu\s u_0)}^2+\nt{\nabla (t_\nu\s v_0)}^2< \sbr{\f{2^*-\p}{2^*-2}\p\nu D_0}^{\f{2}{2-\p}}.
\ee
Since $(u_\nu,v_\nu)\in\PR_\nu^+(a,b)$ we observe that
\be\lab{tem6-14-12}
	\nt{\nabla u_\nu}^2+\nt{\nabla v_\nu}^2< \sbr{\f{2^*-\p}{2^*-2}\p\nu D_0}^{\f{2}{2-\p}}.
\ee
Now using $t_\nu\s(u_0,v_0)$ as a test function of $m_\nu(a,b)$ and by \eqref{tem6-14-13} we get
\be\lab{tem6-14-14}
	\begin{aligned}
	m_\nu(a,b)&\le I_\nu(t_\nu\s(u_0,v_0))=\f{1}{N}(\nt{\nabla (t_\nu\s u_0)}^2+\nt{\nabla (t_\nu\s v_0)}^2)\\
		&\quad\quad	-\nu\f{2^*-\p}{2^*}D_0(\nt{\nabla (t_\nu\s u_0)}^2+\nt{\nabla (t_\nu\s v_0)}^2)^{\f{\p}{2}}.
	\end{aligned}
\ee
Next inequality \eqref{VGNine} implies
\be\lab{tem6-14-15}
	m_\nu(a,b)= I_\nu(u_\nu,v_\nu)
	 \ge\f{1}{N}\big(\nt{\nabla u_\nu}^2+\nt{\nabla v_\nu}^2\big) - \nu\f{2^*-\p}{2^*}D_0\big(\nt{\nabla u_\nu}^2+\nt{\nabla v_\nu}^2\big)^{\f{\p}{2}}.
\ee
A direct calculation shows that the function
	$$f(t)=\f{1}{N}t-\nu\f{2^*-\p}{2^*}D_0 t^{\f{\p}{2}}.$$
is strictly decreasing in $(0,t_0)$, where
	$$t_0:=\sbr{\f{2^*-\p}{2^*-2}\p\nu D_0}^{\f{2}{2-\p}}.$$
Thus \eqref{tem6-14-11}-\eqref{tem6-14-15} yield
\be\lab{est2}
	\nt{\nabla u_\nu}^2+\nt{\nabla v_\nu}^2\ge \nt{\nabla (t_\nu\s u_0)}^2+\nt{\nabla (t_\nu\s v_0)}^2.
\ee
Moreover, similar to \eqref{tem6-14-16}-\eqref{tem6-14-12}, we have
\be\lab{tem6-14-21}
	\nm{|t_\nu\s u_0|^\al |t_\nu\s v_0|^\beta}_1< \sbr{\f{2^*-\p}{2^*-2}\p\nu D_0^{\f{2}{\p}}}^{\f{\p}{2-\p}}
\ee
and
\be
	\nm{u_\nu^\al v_\nu^\beta}_1< \sbr{\f{2^*-\p}{2^*-2}\p\nu D_0^{\f{2}{\p}}}^{\f{\p}{2-\p}}.
\ee
Now using $t_\nu\s(u_0,v_0)$ again  as a test function of $m_\nu(a,b)$ and by \eqref{tem6-14-13},
\be
	\begin{aligned}
	m_\nu(a,b)&\le I_\nu(t_\nu\s(u_0,v_0))
		=\f{1}{N}D_0^{-\f{2}{\p}}\nm{|t_\nu\s u_0|^\al |t_\nu\s v_0|^\beta}_1^{\f{2}{\p}}\\
		&\quad\quad	-\nu\f{2^*-\p}{2^*}\nm{|t_\nu\s u_0|^\al |s_\nu\s v_0|^\beta}_1.
	\end{aligned}
\ee
By inequality \eqref{VGNine} there holds
\be\lab{tem6-14-22}
	m_\nu(a,b) = I_\nu(u_\nu,v_\nu)\ge\f{1}{N}D_0^{-\f{2}{\p}}\nm{u_\nu^\al v_\nu^\beta}_1^{\f{2}{\p}}
		-\nu\f{2^*-\p}{2^*}\nm{u_\nu^\al v_\nu^\beta}_1.
\ee
Next the function
	$$\tilde f(t)=\f{1}{N}D_0^{-\f{2}{\p}}t^{\f{2}{\p}}-\nu\f{2^*-\p}{2^*}t.$$
is strictly decreasing in $(0,\tilde t_0)$, where
	$$\tilde t_0:= \sbr{\f{2^*-\p}{2^*-2}\p\nu D_0^{\f{2}{\p}}}^{\f{\p}{2-\p}}.$$
Therefore \eqref{tem6-14-21}-\eqref{tem6-14-22} give
\be\lab{est3}
	\nm{u_\nu^\al v_\nu^\beta}_1\ge \nm{|t_\nu\s u_0|^\al |t_\nu\s v_0|^\beta}_1.
\ee

	As in Lemma \ref{tem1}, it follows from \eqref{tem6-14-23} and \eqref{tem6-14-24} that
		$$e^{t_\nu}\sim\nu^{\f{1}{2-\p}}\quad\text{as}~\nu\to0.$$
Then we obtain from \eqref{tem6-14-23} that
\be\lab{tem6-14-31}
	e^{t_\nu}=(1+o_\nu(1))\sbr{\f{\nu\p \nm{u_0^\al v_0^\beta}_1}{\nt{\nabla u_0}^2+\nt{\nabla v_0}^2} }^{\f{1}{2-\p}}
		= \big(1+o_\nu(1)\big)\nu^{\f{1}{2-\p}}.
\ee
Since $(u_\nu,v_\nu)$ satisfies \eqref{mainequ}, $(\bar u_\nu,\bar v_\nu) := (-t_\nu)\s(u_\nu,v_\nu)$ satisfies the following equation:
\be\lab{equ1}
	\left\{ 	\begin{aligned}
	&-\dl \bar u_\nu+\la_{1,\nu}e^{-2t_\nu} \bar u_\nu=e^{-(2^*-2)t_\nu}\bar u_\nu^{2^*-1}
		+\al  \nu e^{-(\p-2)t_\nu} \bar u_\nu^{\al-1}\bar v_\nu^\beta&&\quad \text{in }\RN,\\
	&-\dl \bar v_\nu+\la_{2,\nu}e^{-2t_\nu} \bar v_\nu=e^{-(2^*-2)t_\nu}\bar v_\nu^{2^*-1}
		+\beta  \nu e^{-(\p-2)t_\nu} \bar v_\nu^\al \bar v_\nu^{\beta-2}&&\quad \text{in }\RN,\\
	&\int_\RN \bar u_\nu^2=a^2,\int_\RN \bar v_\nu^2=b^2,
	\end{aligned}	\right.
\ee
By Lemma \ref{tem1} we have
	$$\nt{\nabla \bar u_\nu}^2+\nt{\nabla \bar v_\nu}^2=e^{-2t_\nu}(\nt{\nabla u_\nu}^2+\nt{\nabla v_\nu}^2)\sim1.$$
Therefore $\{(\bar u_\nu,\bar v_\nu)\}$ is bounded in $H$.
Recall that $(\bar u_\nu,\bar v_\nu)$ are radial so that we have up to a subsequence,
	$$(\bar u_\nu,\bar v_\nu)\rh(\bar u,\bar v)\quad\text{in} ~H_{rad},$$
	$$(\bar u_\nu,\bar v_\nu)\ra(\bar u,\bar v)\quad\text{in} ~L^q(\RN)\times L^q(\RN),~\text{for}~2<q<2^*,$$
	$$(\bar u_\nu,\bar v_\nu)\ra(\bar u,\bar v)\quad\text{a.e. in}~\RN.$$
Applying Lemma \ref{tem1} once more we see that $\{\la_{1,\nu}e^{-2s_\nu}\}$ and $\{\la_{2,\nu}e^{-2s_\nu}\}$  are bounded, hence
	$$\la_{1,\nu}e^{-2s_\nu}\to\la_1 \quad\text{and}\quad \la_{2,\nu}e^{-2s_\nu}\to\la_2,\quad\text{as}~\nu\to0$$
up to a subsequence. On the other hand, by \eqref{tem6-14-31},
	$$e^{-(2^*-2)t_\nu}\to0,\quad\text{and}\quad \nu e^{-(\p-2)t_\nu}\to1,\quad\text{as}~\nu\to0.$$
Thus, we know that
\be\lab{equ2}
	\left\{ 	\begin{aligned}
	&-\dl \bar u+\la_1 \bar u=\al  \bar u^{\al-1}\bar v^\beta &&\quad \text{in }\RN,\\
	&-\dl \bar v+\la_2 \bar v=\beta \bar v^\al \bar v^{\beta-2} &&\quad \text{in }\RN.
	\end{aligned}	\right.
\ee
Now \eqref{est3} implies $\nm{\bar u^\al \bar v^\beta}_1\ge \nm{u_0^\al v_0^\beta}_1$,
which gives $\bar u\neq0$, $\bar v\neq0$.
By the maximum principle and \cite[Lemma A.2]{Ikoma}, we obtain $\la_1>0$ and $\la_2>0$.
Then combining \eqref{equ1} and \eqref{equ2} we conclude
$(\bar u_\nu,\bar v_\nu)\ra(\bar u,\bar v)$ in $H$ as $\nu\to0$.
It follows that $\nt{\bar u}=a$, $\nt{\bar v}=b$, and by \eqref{est2}:
	$$
		l(a,b)\le K(\bar u,\bar v)=\f{\p-2}{2\p}(\nt{\nabla \bar u}^2+\nt{\nabla \bar v}^2)
			\le \f{\p-2}{2\p}(\nt{\nabla u_0}^2+\nt{\nabla v_0}^2)=K(u_0,v_0)
			=l(a,b).
	$$
We finally deduce $(\bar u,\bar v)\in\GR(a,b)$, i.e.,
	$$dist_H\sbr{(-t_\nu)\s(u_\nu,v_\nu),\GR(a,b)}\to0,\quad\text{as}~\nu\to0^+.$$
\ep

	Now we turn to the case $\al+\beta\ge2+\f{4}{N}$.
\bl\lab{tem5}
	If $\al+\beta\ge 2+\f{4}{N}$ then
		$$dist_{\DR^{1,2}(\RN)}\big((u_\nu,v_\nu),\UR\times\{0\}\big) \to 0 \quad\text{as}~\nu\to0^+,$$
	or
		$$dist_{\DR^{1,2}(\RN)}\big((u_\nu,v_\nu),\{0\}\times\UR\big) \to 0 \quad\text{as}~\nu\to0^+;$$
	here $\UR=\lbr{U_\e:\e>0}$.
\el

\bp
	First observe that$\{(u_\nu,v_\nu)\}$ is bounded in $H$ because $(u_\nu,v_\nu)\in\PR_\nu(a,b)$ satisfies $I_\nu(u_\nu,v_\nu)=m_\nu(a,b)$.
Let $\nt{\nabla u_\nu}\to l_1$ and $\nt{\nabla v_\nu}\to l_2$ with $l_1\ge0$, $l_2\ge0$.
If $l_1=l_2=0$ then
	$$0=\lim_{\nu\to0}I_\nu(u_\nu,v_\nu)=\lim_{\nu\to0}m_\nu(a,b)>0,$$using 
which is a contradiction. The Sobolev inequality and $P_\nu(u_\nu,v_\nu)\to0$ imply
	$$l_1^2+l_2^2\le \SR^{-2^*/2}\big(l_1^{2^*}+l_2^{2^*}\big),$$
hence $l_1^2+l_2^2\ge \SR^{N/2}$ and equality holds if and only if $l_1=0$ or $l_2=0$.
Consequently we have
	$$\f{1}{N}\SR^{N/2}\ge \lim_{\nu\to0}I_\nu(u_\nu,v_\nu)=\f{1}{N}(l_1^2+l_2^2)$$
and therefore either $l_1=0$, $l_2=\SR^{N/4}$, or $l_1=\SR^{N/4}$, $l_2=0$.
If $l_1=0$, $l_2=\SR^{N/4}$ then
	$$\nt{\nabla u_\nu}^2\to0,\quad \nt{\nabla v_\nu}^2\to\SR^{N/2},\quad
		\ns{v_\nu}^{2^*}\to\SR^{N/2} \quad\text{as}~\nu\to0,$$
hence $v_\nu$ is a minimizing sequence of \eqref{Sobolev1}. Now \cite[Theorem 1.41]{Willem} implies that there exists $\q_\nu>0$ such that for some $\e_*>0$
	$$\Big(u_\nu,\q_\nu^{\f{N-2}{2}}v_\nu(\q_\nu\cdot)\Big)\to(0,U_{\e_*})\quad\text{strongly in}~\DR^{1,2}(\RN)~\text{as}~\nu\to0^+,$$
up to a subsequence. The case $l_1=\SR^{N/4}$, $l_2=0$ can be treated analogously.
\ep

\bp[Proof of Theorem \ref{thm3}]
	This follows immediately from Lemma \ref{tem1}, Lemma \ref{tem4} and Lemma \ref{tem5}.
\ep

\bp[Proof of Theorem \ref{thm4}]
	Since the proof is similar to that of Lemma \ref{tem4}, we only sketch it.
Given $(u,v)\in\TR(a,b)$ we use $t_\nu(u,v)\s(u,v)$ as a test function for $m_\nu(a,b)$
where $t_\nu(u,v)$ is given by Lemma \ref{structure2}.
Using similar arguments as in the proof of Lemma \ref{tem1} and direct calculations we get
$m_\nu(a,b)\lesssim \nu^{\f{2}{2-\p}}$ as $\nu\to+\iy$.
It follows that
	$$\nt{\nabla u_\nu}^2+\nt{\nabla v_\nu}^2 \to 0 \quad\text{and}\quad \nm{u_\nu^\al v_\nu^\beta}_1 \to 0\quad\text{as}~\nu\to+\iy.$$
Since
\be
	\nt{\nabla u_\nu}^2+\nt{\nabla v_\nu}^2=\ns{u_\nu}^{2^*}+\ns{v_\nu}^{2^*}+\p\nu\nm{u_\nu^\al v_\nu^\beta}_1,
\ee
using the Sobolev inequality and \eqref{VGNine}, we obtain 
	$$\nt{\nabla u_\nu}^2+\nt{\nabla v_\nu}^2\ge \big(1+o_\nu(1)\big)\sbr{\nu \p D_0}^{\f{2}{2-\p}}$$
and
	$$\nm{u_\nu^\al v_\nu^\beta}_1\ge \big(1+o_\nu(1)\big)\sbr{\nu \p D_0^{\f{2}{\p}}}^{\f{\p}{2-\p}},$$
where $D_0$ is defined by \eqref{est1}.
For any $(u_0,v_0)\in\GR(a,b)$, let $t_\nu:=t_\nu(a,b)$ be given by Lemma \ref{structure2}.
Then
	$$e^{2t_\nu}\big(\nt{\nabla u_0}^2+\nt{\nabla v_0}^2\big)=e^{2^*t_\nu}\big(\ns{u_0}^{2^*}+\ns{v_0}^{2^*}\big)
			+e^{\p t_\nu}\p\nu\nm{|u_0|^\al |v_0|^\beta}_1.$$
It follows from \eqref{est1} that
	$$e^{2t_\nu}\big(\nt{\nabla u_0}^2+\nt{\nabla v_0}^2\big)\le\sbr{\nu \p D_0}^{\f{2}{2-\p}},$$
and
	$$e^{\p t_\nu}\nm{u_0^\al v_0^\beta}_1\le  \sbr{\nu \p D_0^{\f{2}{\p}}}^{\f{\p}{2-\p}}.$$
Thus we get
\be
	\begin{aligned}
	I_\nu\big(t_\nu\s(u_0,v_0)\big)&=\f{\p-2}{2\p}e^{2t_\nu}\big(\nt{\nabla u_0}^2+\nt{\nabla v_0}^2\big)
		         + \f{2^*-\p}{2^*\p}e^{2^*t_\nu}\big(\ns{u_0}^{2^*}+\ns{v_0}^{2^*}\big)\\
		&\le \big(1+o_\nu(1)\big)\f{\p-2}{2\p}\sbr{\nu \p D_0}^{\f{2}{2-\p}}.
	\end{aligned}
\ee
Note that $I_\nu(t_\nu\s(u_0,v_0)\ge m_\nu(a,b)$ and
	$$m_\nu(a,b)=I_\nu(u_\nu,v_\nu) = \big(1+o_\nu(1)\big)\f{\p-2}{2\p}\sbr{\nu \p D_0}^{\f{2}{2-\p}},$$
as $\nu\to+\iy$, we must have
\be\lab{tem6-15-1}
	\nt{\nabla u_\nu}^2+\nt{\nabla v_\nu}^2 = \big(1+o_\nu(1)\big)\sbr{\nu \p D_0}^{\f{2}{2-\p}}.
\ee
Similarly, we also have
	$$\nm{u_\nu^\al v_\nu^\beta}_1=(1+o_\nu(1))\sbr{\nu \p D_0^{\f{2}{\p}}}^{\f{\p}{2-\p}}.$$
Moreover, using $K(u_0,v_0)=l(a,b)$ and $L(u_0,v_0)=0$, we have
	$$\nm{u_0^\al v_0^\beta}_1=\sbr{\p D_0^{\f{2}{\p}}}^{\f{\p}{2-\p}},$$
which implies
\be\lab{tem6-15-2}
	\nm{u_\nu^\al v_\nu^\beta}_1=(1+o_\nu(1))\nu^{\f{\p}{2-\p}}\nm{u_0^\al v_0^\beta}_1.
\ee
Let $(\bar u_\nu,\bar v_\nu)=s_\nu\s(u_\nu,v_\nu)$ with $e^{s_\nu}:=\nu^{\f{1}{\p-2}}$.
Then
	$\nt{\nabla \bar u_\nu}^2=a^2,\quad \nt{\nabla \bar v_\nu}^2=b^2,$
and from \eqref{tem6-15-1},
	$$\nt{\nabla \bar u_\nu}^2+\nt{\nabla \bar v_\nu}^2 = \big(1+o_\nu(1)\big)\sbr{ \p D_0}^{\f{2}{2-\p}}.$$
Therefore $\{(\bar u_\nu,\bar v_\nu)\}$ is bounded in $H$. Since $(\bar u_\nu,\bar v_\nu)$ are radial, we have up to a subsequence:
	$$(\bar u_\nu,\bar v_\nu)\rh(\bar u,\bar v)\quad\text{in} ~H_{rad},$$
	$$(\bar u_\nu,\bar v_\nu)\ra(\bar u,\bar v)\quad\text{in} ~L^q(\RN)\times L^q(\RN),~\text{for}~2<q<2^*,$$
	$$(\bar u_\nu,\bar v_\nu)\ra(\bar u,\bar v)\quad\text{a.e. in}~\RN.$$
Moreover,
\be
	\left\{ 	\begin{aligned}
	&-\dl \bar u_\nu+\la_{1,\nu}e^{2s_\nu} \bar u_\nu=e^{-(2^*-2)s_\nu}\bar u_\nu^{2^*-1}
		+\al    \bar u_\nu^{\al-1}\bar v_\nu^\beta,\quad \text{in }\RN,\\
	&-\dl \bar v_\nu+\la_{2,\nu}e^{2s_\nu} \bar v_\nu=e^{-(2^*-2)s_\nu}\bar v_\nu^{2^*-1}
		+\beta  \bar v_\nu^\al \bar v_\nu^{\beta-2},\quad \text{in }\RN,\\
	&\int_\RN \bar u_\nu^2=a^2,\int_\RN \bar v_\nu^2=b^2.
	\end{aligned}	\right.
\ee
Since
	$$\la_{1,\nu}a^2+\la_{2,\nu}b^2=(1-\p)\nu \nm{u_\nu^\al v_\nu^\beta}_1\sim \nu^{\f{2}{\p-2}},$$
 by the Pohozaev identity, we see that $\la_{1,\nu}e^{2s_\nu}$ and $\la_{2,\nu}e^{2s_\nu}$ are bounded.
Then as in Lemma \ref{tem4}, we obtain that
	$$dist_H\sbr{s_\nu\s(u_\nu,v_\nu),\GR(a,b)}\to0,$$
as $\nu\to+\iy$.
\ep

\vskip0.4in
\begin{appendices}

\section{Vector-valued Gagliardo-Nirenberg inequality}

Recall the Gagliardo-Nirenberg inequality,
\be\lab{GNine}
	\nm{u}_p^p\le \CR(N,p)\nt{u}^{p-\ga_p}\nt{\nabla u}^{\ga_p}\quad\text{for all $u\in H^1(\RN)$,}
\ee
where the sharp constant $\CR(N,p)$ satisfies
\be\lab{GNconstant}
	\CR(N,p)^{-1}=\inf_{u\in H^1(\RN)\setminus\{0\}} \f{\nt{u}^{p-\ga_p}\nt{\nabla u}^{\ga_p}}{\nm{u}_p^p}
			=\ga_p^{\f{\ga_p}{2}}(1-\ga_p)^{1-\f{\ga_p}{2}}\nt{Z}^{p-2}.
\ee
Here $Z$ is the unique solution of
\be\lab{Equ}
	\left\{ \begin{aligned}
	&-\Delta w+w=|w|^{p-2}w &&\quad\text{in}~\RN,\\
	& w>0\ \text{and}\ w(x)\to0 &&\quad\text{as}~|x|\to\iy,\\
	& w(0)=\max_{x\in\RN}w(x).
	\end{aligned} \right.
\ee
For a detailed proof we refer to \cite{Weinstein=1982}. Moreover, the function $Z^{\sigma,\mu}(x):=\sigma Z(\mu x)$ satisfies
	$$-\Delta Z^{\sigma,\mu}+\mu^2Z^{\sigma,\mu}=\sigma^{2-p}\mu^2|Z^{\sigma,\mu}|^{p-2}Z^{\sigma,\mu},$$
and \eqref{GNconstant} is achieved by $u$ if and only if
	$$u\in\lbr{ Z^{\sigma,\mu}(\cdot+y): \sigma>0, \mu>0, y\in\RN }.$$
When considering systems we need a vector-valued version of the Gagliardo-Nirenberg inequality.
For $\al>1$, $\beta>1$, $\al+\beta<2^*$ we set
	$$Q(u,v) := \f{ \sbr{\nt{u}^2+\nt{v}^2}^{(\al+\beta-\p)/2}
		\sbr{\nt{\nabla u}^2+\nt{\nabla v}^2}^{\p/2}} {\nm{|u|^\al|v|^\beta}_1}.$$
and
\be\lab{VGNconstant}
	\CR(N,\al,\beta)^{-1}:=   \inf_{u,v\in H^1(\RN)\setminus\{0\}} Q(u,v).
\ee
It follows easily from \eqref{GNine} that $0<\CR(N,\al,\beta)^{-1}<+\iy$. From definition \eqref{VGNconstant} we obtain the vector-valued Gagliardo-Nirenberg inequality
\be\lab{VGNine}
	\big\||u|^\al|v|^\beta\big\|_1\le \CR(N,\al,\beta)\sbr{\nt{u}^2+\nt{v}^2}^{(\al+\beta-\p)/2}
			 \sbr{\nt{\nabla u}^2+\nt{\nabla v}^2}^{\p/2},
\ee
which holds for $u,v\in H^1(\RN)$. The vector-valued Gagliardo-Nirenberg inequality has been studied by many researchers,
see \cite{VGN-1,VGN-2,VGN-3} and the reference therein.

\begin{proposition}\lab{VGN}
	The sharp constant $\CR(N,\al,\beta)$ is determined by
		$$\CR(N,\al,\beta)=\mbr{ \sbr{\f{\al}{\beta}}^{\f{\beta}{\al+\beta}}
				+ \sbr{\f{\beta}{\al}}^{\f{\al}{\al+\beta}} }^{-\f{\al+\beta}{2}} \CR(N,\al+\beta).$$
	Moreover, $\CR(N,\al,\beta)^{-1}$ is achieved by $(u,v)$ if and only if
	\be\lab{A1}
		(u,v)\in\lbr{ \big(Z_1^{\sigma,\mu}(\cdot+y), Z_2^{\sigma,\mu}(\cdot+y)\big): \sigma>0, \mu>0, y\in\RN },
	\ee
	where
		$$(Z_1,Z_2)=\sbr{\al^{-\f{2-\beta}{2(\al+\beta-2)}}\beta^{-\f{\beta}{2(\al+\beta-2)}},
				 \al^{-\f{\al}{2(\al+\beta-2)}}\beta^{-\f{2-\al}{2(\al+\beta-2)}}} Z(x)$$
	is the unique solution of
	\be
		\left\{  \begin{aligned}
		&-\Delta u+u=\al u^{\al-1}v^\beta &&\quad\text{in}~\RN,\\
		&-\Delta v+v=\beta u^\al v^{\beta-1} &&\quad\text{in}~\RN,\\
		&u>0, v>0,\quad u,v\in H^1(\RN).
		\end{aligned}  \right.
	\ee
\end{proposition}
\bp
Let $(u_n,v_n)$ be a minimizing sequence of $\CR(N,\al,\beta)^{-1}$, i.e., $Q(u_n,v_n)\to\CR(N,\al,\beta)^{-1}$ as $n\to\iy$.
It is easy to check that $Q(u^{\sigma,\mu},v^{\sigma,\mu})=Q(u,v)$ for any $\sigma,\mu>0$.
Then by scaling and rearrangement, we may assume that
$u_n,v_n\in H_{rad}^1(\RN)$, $u_n,v_n\ge0$ and
	$$\nt{\nabla u_n}^2+\nt{\nabla v_n}^2=\nt{u_n}^2+\nt{v_n}^2=1,\quad\text{for all}~n\ge1.$$
Passing to a subsequence we get $(u_n,v_n)\rh(u,v)$ in $H_{rad}$. Since $2<\al+\beta<2^*$, we have
	$$ \nm{u_n^\al v_n^\beta}_1\to \nm{u^\al v^\beta}_1,$$
which means that $\nm{u^\al v^\beta}_1=\CR(N,\al,\beta)$ and $u\neq0$, $v\neq0$.
Noting that
	$$\nt{\nabla u}^2+\nt{\nabla v}^2\le \lim_{n\to\iy} \big(\nt{\nabla u_n}^2+\nt{\nabla v_n}^2\big) = 1$$
and
	$$\nt{u}^2+\nt{v}^2\le \lim_{n\to\iy} \big(\nt{u_n}^2+\nt{v_n}^2\big) = 1,$$
we obtain
	$$\CR(N,\al,\beta)^{-1}\le Q(u,v)\le \nm{u^\al v^\beta}_1^{-1}=\CR(N,\al,\beta)^{-1}.$$
It follows that $(u_n,v_n)\ra(u,v)$ in $H$ and $Q(u,v)=\CR(N,\al,\beta)^{-1}$.
Moreover, by $Q'(u,v)=0$ we conclude that
\be
	\begin{cases}
	-\Delta u+\f{\al+\beta-\p}{\p}u=\f{\al \CR(N,\al,\beta)^{-1}}{\p}u^{\al-1}v^\beta \quad&\text{in}~\RN,\\
	-\Delta v+\f{\al+\beta-\p}{\p}v=\f{\beta \CR(N,\al,\beta)^{-1}}{\p}u^\al v^{\beta-1} \quad&\text{in}~\RN.
	\end{cases}
\ee
Hence by maximum principle $u>0$, $v>0$.
By scaling, there exists $\sigma_0,\mu_0>0$ such that
\be
	\begin{cases}
	-\Delta u^{\sigma_0,\mu_0}+ u^{\sigma_0,\mu_0}=
			\al (u^{\sigma_0,\mu_0})^{\al-1}(v^{\sigma_0,\mu_0})^\beta, \quad&\text{in}~\RN,\\
	-\Delta v^{\sigma_0,\mu_0}+ v^{\sigma_0,\mu_0}=
			\beta (u^{\sigma_0,\mu_0})^\al (v^{\sigma_0,\mu_0})^{\beta-1},\quad&\text{in}~\RN.
	\end{cases}
\ee
Setting
	$$u^{\sigma_0,\mu_0}=\al^{-\f{2-\beta}{2(\al+\beta-2)}}\beta^{-\f{\beta}{2(\al+\beta-2)}} U(x),$$
and
	$$v^{\sigma_0,\mu_0}=\al^{-\f{\al}{2(\al+\beta-2)}}\beta^{-\f{2-\al}{2(\al+\beta-2)}} V(x),$$
we have
$$\begin{cases}
	-\Delta U+ U=U^{\al-1}V^\beta \quad&\text{in}~\RN,\\
	-\Delta V+ V=U^\al V^{\beta-1}\quad&\text{in}~\RN.
\end{cases}$$
Following the proof of \cite[Theorem 1]{VGN-1}, we can prove $U=V=Z(x)$.
Thus we have
	$$u^{\sigma_0,\mu_0}=\al^{-\f{2-\beta}{2(\al+\beta-2)}}\beta^{-\f{\beta}{2(\al+\beta-2)}} Z(x),$$
	$$v^{\sigma_0,\mu_0}=\al^{-\f{\al}{2(\al+\beta-2)}}\beta^{-\f{2-\al}{2(\al+\beta-2)}} Z(x),$$
and
	$$\CR(N,\al,\beta)^{-1}=Q(u,v)=Q(u^{\sigma_0,\mu_0},v^{\sigma_0,\mu_0})=
		\mbr{ \sbr{\f{\al}{\beta}}^{\f{\beta}{\al+\beta}}
				+ \sbr{\f{\beta}{\al}}^{\f{\al}{\al+\beta}} }^{\f{\al+\beta}{2}} \CR(N,\al+\beta)^{-1}.$$
Finally, for any minimizer $(u,v)$ of $\CR(N,\al,\beta)^{-1}$, repeating the above process we can prove \eqref{A1}.
\ep

\end{appendices}


{\sc Address of the authors:}\\

\noindent
Thomas Bartsch\\
Mathematisches Institut \\
University of Giessen \\
Arndtstr.\ 2 \\
35392 Giessen \\
Germany \\
Thomas.Bartsch@math.uni-giessen.de \\

\noindent
Houwang Li \& Wenming Zou\\
Department of Mathematical Sciences\\
Tsinghua University\\
Beijing 100084\\
China\\
li-hw17@mails.tsinghua.edu.cn\quad\&\quad zou-wm@mail.tsinghua.edu.cn\\


 \end{document}